\documentstyle[12pt]{article}
\parskip=\medskipamount                 
\thicklines                     
\newcommand{\bim}[6]{\bibitem{#1}#2, {\em #3\/}$\;${\bf
#4}$\;$(#5)$\;${#6}.}
\def\IR{\relax{\rm I\kern-.18em R}}
\def\ZZ{\relax{\sf Z\kern-.4em Z}}
\def\a{\alpha} \def\b{\beta}    
  \def\om{\omega}
   
   \def\cD{{\cal D}}
   
 \def\cK{{\cal K}} \def\cL{{\cal L}}

\newtheorem{proposition}{Proposition}[section]
\newtheorem{corollary}{Corollary}[section]

\newtheorem{lemma}{Lemma}[section]
\newtheorem{definition}{Definition}[section]

\catcode`\@=11

\newif\if@fewtab\@fewtabtrue

\catcode`\@=11

\newif\if@fewtab\@fewtabtrue

{\count255=\time\divide\count255 by 60
\xdef\hourmin{\number\count255}
\multiply\count255 by-60\advance\count255 by\time
\xdef\hourmin{\hourmin:\ifnum\count255<10 0\fi\the\count255}}
\def\ps@draft{\let\@mkboth\@gobbletwo
    \def\@oddhead{}
    \def\@oddfoot
       {\hbox to 7 cm{$\scriptstyle Draft\ version:\ \draftdate$
       \hfil}\hskip -7cm\hfil\rm\thepage \hfil}
    \def\@evenhead{}\let\@evenfoot\@oddfoot}

\def\ceqno{\global\@fewtabfalse
    \ifcase\@eqcnt \def\@tempa{& & &}\or \def\@tempa{& &}
      \or \def\@tempa{&}
      \or\def\@tempa{}\fi\@tempa
{\rm(\theequation)}}

\def\aeqno#1{\global\@fewtabfalse
    \ifcase\@eqcnt \def\@tempa{& & &}\or \def\@tempa{& &}
      \or \def\@tempa{&}
      \or\def\@tempa{}\fi\@tempa
{\rm(\theequation,#1)}}

\def\label#1{\ifnum\draftcontrol=1
 \global\def\draftnote{$\scriptstyle #1$}\fi
 \@bsphack\if@filesw {\let\thepage\relax
   \def\protect{\noexpand\noexpand\noexpand}%
\xdef\@gtempa{\write\@auxout{\string
      \newlabel{#1}{{\@currentlabel}{\thepage}}}}}\@gtempa
   \if@nobreak \ifvmode\nobreak\fi\fi\fi
  \@esphack}

\def\alabel#1#2{\label{#1}\global\@fewtabfalse
    \ifcase\@eqcnt \def\@tempa{& & &}\or \def\@tempa{& &}
      \or \def\@tempa{&}
      \or\def\@tempa{}\fi\@tempa
{\hbox to 3cm{\phantom{\rm(\theequation,#2)}
\draftnote \hfil}\hskip -3cm {\rm(\theequation,#2)}}}

\def\clabel#1{\label{#1}\global\@fewtabfalse
    \ifcase\@eqcnt \def\@tempa{& & &}\or \def\@tempa{& &}
      \or \def\@tempa{&}
      \or\def\@tempa{}\fi\@tempa
{\hbox to 3cm{\phantom{\rm(\theequation)}
\draftnote \hfil}\hskip -3cm{\rm(\theequation)}}}

\def\eqnarray{\def\draftnote{{}}\global\@fewtabtrue
\stepcounter{equation}\let\@currentlabel=\theequation
\global\@eqnswtrue
\global\@eqcnt\z@\tabskip\@centering\let\\=\@eqncr
$$\halign to \displaywidth\bgroup\@eqnsel\hskip\@centering\@eqcnt\z@
  $\displaystyle\tabskip\z@{##}$&\global\@eqcnt\@ne
  \hskip 1\arraycolsep \hfil${##}$\hfil
  &\global\@eqcnt\tw@ \hskip 1\arraycolsep
$\displaystyle\tabskip\z@{##}$
\hfil  \tabskip\@centering&\global\@eqcnt\thr@@\llap{##}\tabskip\z@
\cr}

\def\endeqnarray{\@@eqncr\egroup
      \global\advance\c@equation\m@ne$$\global\@ignoretrue}

\def\@eqnnum{\hbox to 3cm{\phantom{\rm(\theequation)} \draftnote
                         \hfil}\hskip -3cm {\rm(\theequation)}}

\def\@@eqncr{\let\@tempa\relax
    \ifcase\@eqcnt \def\@tempa{& & &}\or \def\@tempa{& &}
      \or \def\@tempa{&}
      \or\def\@tempa{}
\fi\@tempa
\if@eqnsw
\if@fewtab\@eqnnum\fi
\stepcounter{equation}\fi\global
\@eqnswtrue\global\@eqcnt\z@\global\@fewtabtrue\cr}

\def\draftcite#1{\ifnum\draftcontrol=1#1\else{}\fi}

\def\@lbibitem[#1]#2{\item{}\hskip -3cm \hbox to 2cm
{\hfil$\scriptstyle\draftcite{#2}$}\hskip
1cm[\@biblabel{#1}]\if@filesw
     {\def\protect##1{\string ##1\space}\immediate
      \write\@auxout{\string\bibcite{#2}{#1}}}\fi\ignorespaces}

\def\@bibitem#1{\item\hskip -3cm \hbox to 2cm
{\hfil $\scriptstyle\draftcite{#1}$}\hskip 1cm
\if@filesw \immediate\write\@auxout
       {\string\bibcite{#1}{\the\value{\@listctr}}}\fi\ignorespaces}

\def\nsection#1{\section{#1}\setcounter{equation}{0}}

\def\draftdate{\number\month/\number\day/\number\year\ \ \ \hourmin }

\global\def\draftcontrol{0}
\catcode`\@=12

\def\theequation{{\thesection.\arabic{equation}}}

\def\qq{\begin{eqnarray}}
\def\qqq{\end{eqnarray}}

\newlength{\shiftwidth}
\def\shift#1{&&\hbox to \shiftwidth{\hfill $\displaystyle#1$}}
\newlength{\sshiftwidth}
\def\sshift#1{\lefteqn{\hbox to
\sshiftwidth{\hfill$\displaystyle#1$}}}
\def\llefteqn#1{\hbox to 0pt{$\displaystyle #1 $\hss}\hspace*{1in}}

\def\ie{{\it i.e.\ }}
\def\eg{{\it e.g.\ }}
\def\cf{{\it cf.\ }}
\def\rhs{{\it r.h.s.\ }}

\def\Rhs{RHS\ }
\def\ihs{$\ZZ$HS\ }
\def\asl{ASL\ }
\def\sasl{SASL\ }
\def\bl{BL\ }

\def\Tr{\mathop{{\rm Tr}}\nolimits}

\def\ord{\mathop{{\rm ord}}\nolimits}
\def\p{^{\prime}}
\def\pp{^{\prime\prime}}
\def\Tub{\mathop{{\rm Tub}}\nolimits}
\def\sign#1{\mathop{{\rm sign}\left(#1\right)}\nolimits}
\def\si{\mathop{{\rm si}}\nolimits}
\def\spint#1{\int\limits^{+\infty}_{\scriptstyle -\infty \atop [{#1} =
0]}}

\def\lk{\mathop{{\rm lk}}\nolimits}
\def\max{\mathop{{\rm max}}\nolimits}
\def\mmax#1{\max\{#1\}}

\def\Pexp{\mathop{{\rm Pexp}}\nolimits}
\def\PexpA#1{\Pexp \left(\oint_{#1} A_\mu dx^\mu \right)}
\def\cs{S_{\rm CS}}
\def\csc{S^{(c)}_{\rm CS}}
\def\ordH{\ord H_1(M,\ZZ)}
\def\ordHp{\ord H_1(M\p,\ZZ)}

\def\aa{{\a_1,\ldots,\a_N}}
\def\daa{d\a_1 \cdots d\a_N}
\def\da{da_1 \cdots da_N}
\def\paa{(a_1,\ldots,a_N)}

\def\ztr{Z^{\rm (tr)}}
\def\zaml{Z_\aa(M,\cL;k)}
\def\zasl{Z_\aa(S^3,\cL;k)}
\def\zacl{Z^{(c)}_\aa(M,\cL;k)}
\def\ztrm{Z^{\rm (tr)}(M;k)}
\def\ztram{Z^{\rm (tr)}_\aa(M;k)}
\def\ztras{Z^{\rm (tr)}_\aa(S^3;k)}
\def\jas{J_\aa(S^3,\cL;k)}
\def\zmpk{Z(M\p;k)}
\def\ztrmpk{Z^{\rm (tr)}(M\p;k)}

\def\jtram{J^{\rm (tr)}_\aa(M,\cL;k)}
\def\jtrass{J^{\rm (tr)}_{2,\ldots 2}(S^3,\cL;k)}
\def\jttram{\tilde{J}^{\rm (tr)}_\aa(M,\cL;k)}

\def\jtrams{J^{\rm (tr)}_\a(M,\cL;k)}
\def\jtraks{J^{\rm (tr)}_{Ka}(M,\cL;k)}
\def\jttraks{\tilde{J}^{\rm (tr)}_{Ka}(M,\cL;k)}
\def\jttrakN{\tilde{J}^{\rm (tr)}_{Ka_1,\ldots,Ka_N}(M,\cL;k)}

\def\kcomp{M\setminus \Tub{\cK}}

\def\Am{\Delta_A \left( M,\cK; e^{2\pi i a} \right)}
\def\Amd{\Delta_A \left( M,\cK; e^{2\pi i {a\over m_2 d} } \right)}
\def\AtN{\tilde{\Delta}_A
    \left( M, \cL; e^{2\pi i a_1}, \ldots, e^{2\pi i a_N} \right)}

\def\RN{\tau_R \left(M \setminus \Tub{\cL}; e^{2\pi i
	a_1},\ldots,e^{2\pi i a_N} \right) }
\def\AtmdN{\tilde{\Delta}_A
    \left( M, \cL; e^{2\pi i {a_1 \over m_2^{(1)} d^{(1)}}},
       \ldots, e^{2\pi i {a_N \over m_2^{(N)} d^{(N)}}}
		\right) }

\def\Rcomp{\tau_R \left( \kcomp; e^{2\pi i a} \right)}

\def\pjN{\prod_{j=1}^N}
\def\sjN{\sum_{j=1}^N}

\def\smz{\sum_{m\geq 0}}
\def\snf{\sum_{-{2\over 3}m \leq n < \infty}}

\def\snzi{\sum_{n=0}^\infty}
\def\smzi{\sum_{m=0}^\infty}
\def\snoi{\sum_{n=1}^\infty}
\def\smtw{\sum_{m=2}^\infty}
\def\smth{\sum_{m=3}^\infty}
\def\slmp{\sum_{l,m = 0 \atop l + m \neq 0}^\infty}

\def\ppq{{(p,q)}}
\def\ppqj{{(p_j,q_j)}}
\def\upq{U^\ppq}
\def\uhpq{\hat{U}^\ppq}
\def\utpq{\tilde{U}^\ppq}
\def\upqj{U^\ppqj}
\def\uhpqj{\hat{U}^\ppqj}
\def\utpqj{\tilde{U}^\ppqj}

\def\snm{S_n(M)}
\def\snmp{S_n(M\p)}
\def\stnlm{\tilde{S}_n(\cL,M)}

\def\sint{S^{\rm (int)}_n(M)}

\def\sns{\sum_{n=1}^\infty \snm \ipK^n}
\def\som{S_1(M)}
\def\lcw{\lambda_{\rm CW}}
\def\opr{${\mbox{Ohtsuki}}^\prime\;\;$}

\def\ffr{\phi_{\rm fr}}
\def\Dfr{\Delta_{\rm fr}}

\def\Dn{\Delta_n}

\def\dt{\tilde{d}}
\def\Dt{\tilde{D}}
\def\dmn{d_{m,n}}
\def\dtmn{\dt_{m,n}}
\def\dpmn{d^{(m,n)}_{m_1,\ldots,m_N}}
\def\dtpmn{\dt^{(m,n)}_{m_1,\ldots,m_N}}
\def\Dmn{D_{m,n}}
\def\Dpmn{D^{(m,n)}_{m_1,\ldots,m_n}}
\def\Dtmn{\Dt_{m,n}}
\def\Dtpmn{\Dt^{(m,n)}_{m_1,\ldots,m_n}}
\def\cnm{C^{(n)}_{m_1,\ldots,m_N}}
\def\cnmp{C^{(n)}_{m_1\p,\ldots,m_N\p}}
\def\cmms{C^{(n+\sjN m_j)}_{m_1,\ldots,m_N}}

\def\Dspmn{D^{(m,2m+n)}_{m_1,\ldots,m_n}}
\def\dz{d^{(m,n)}_{0,m_2,\ldots,m_N}}

\def\pqj{{p_j\over q_j}}
\def\pqlj{p_j + q_j l_{jj}}
\def\ffpqlj{ {p_j \over q_j} + l_{jj} }
\def\fpqlj{ \left( {p_j \over q_j} + l_{jj} \right) }
\def\pqlz{p + q l_{00}}
\def\ffpqlz{ {p \over q} + l_{00} }
\def\fpqlz{ \left( {p \over q} + l_{00} \right) }
\def\prql{\pjN(\pqlj)}

\def\sjn{ \sum_{j=1}^n }
\def\pjn{ \prod_{j=1}^n }
\def\snchi{ S_n(\chi_\cK(S^3) ) }
\def\qmu{ q + \sjn \mu_j }
\def\pqmu{ \left( p, q + \sjn\mu_j \right) }
\def\smup{ \sum_{\mu_j = \pm 1 \atop (1\leq j\leq n) }
    \left( \pjn \mu_j \right) }
\def\compchi{ \chi_{\cK,\pqmu} }

\def\dgr{{\rm dgr}}
\def\Dgrn{{\rm Dgr}_n}

\def\Dgrnl{{\rm Dgr}_n(\cL)}
\def\Dgrnol{{\rm Dgr}_{n+1}(\cL)}
\def\Ddgr{\Delta_\dgr}

\def\va{\vec{a}}
\def\vx{\vec{x}}
\def\vn{\vec{n}}
\def\vs{\vec{\sigma}}

\def\pvaj{\pjN \left( {K\over 4\pi} {d^2 \va_j \over |\va_j|}
        \right)}
\def\pjq{\left( \pjN {\sign{q_j}\over \sqrt{|q_j|}} \right)}
\def\pvtj{\left( \pjN {d^3\va_j\over 4\pi} \right) }
\def\ida{{d^2\va\over 4\pi|\va|}}

\def\pjqs{\frac {\sin \left( \pi {|\va_j| \over q_j} \right)}
		     {|\va_j|}}
\def\pva{(\va_1,\ldots,\va_N)}
\def\lm{L_m\pva}
\def\pml{P_{m,l}\pva}
\def\pe{\left[ 1 + \slmp K^{-l} \pml \right]}

\def\Lt{L^{\rm (tot)}}
\def\lij{l_{ij}}
\def\ljj{l_{jj}}
\def\lmu{l^{(\mu)}}
\def\Mt{\lmu_{ijk}}
\def\Mq{\lmu_{ijkl}}
\def\Mj{\lmu_{j_1,\ldots,j_m}}
\def\Mqij{\lmu_{iijj}}
\def\Mts{\left( \Mt \right)^2}

\def\mij{M_{ij,\mu\nu}}
\def\ppjj{\left(p_{jj} + {\pi^2 \over 6} \right)}

\def\om{\Omega_{1,1}}
\def\oxy#1{\om(x_1,y_{#1})\om(x_2,y\p_{#1})}
\def\oxx{\om(x_1,x_2)}
\def\oiky#1{\oint_{\cK_{#1}} dy_{#1}}
\def\oikyp#1{\oint_{\cK_{#1}} dy\p_{#1}}

\def\gev{G_{\rm ev}}
\def\ipik{{i\pi\over 2K}}
\def\ipiK{{i\pi K\over 2}}
\def\ipK{ \left( {i\pi \over K} \right) }
\def\Kmn{K^{-n}}
\def\sil{\si(\cL)}
\def\ss{{n\over 1-\sil}}
\def\lmn{m+l, n-2l}
\def\eipk{e^{i\pi\over K}}
\def\emipk{e^{-{i\pi\over K}}}

\begin{document}

\begin{titlepage}
\centerline{\hfill                 UMTG-182-95}
\centerline{\hfill                 q-alg/9503011}
\vfill
\begin{center}

{\large \bf
The Trivial Connection Contribution to Witten's Invariant and Finite
Type Invariants of Rational Homology Spheres.
} \\

\bigskip
\centerline{L. Rozansky\footnote{Work supported
by the National Science Foundation
under Grant No. PHY-92 09978.
}}

\centerline{\em Physics Department, University of Miami
}
\centerline{\em P. O. Box 248046, Coral Gables, FL 33124, U.S.A.}
\centerline{{\em E-mail address: rozansky@phyvax.ir.miami.edu}}

\vfill
{\bf Abstract}

\end{center}
\begin{quotation}

We derive an analog of Melvin-Morton bound on the power series
expansion of Jones polynomial of algebraically split links and
boundary links. This allows us to produce a simple formula for the
trivial connection contribution to Witten's invariant of rational
homology spheres. We show that the $n$th term in the $1/K$ expansion
of the logarithm of this contribution is a finite type invariant of
Ohtsuki order $3n$ and of at most Garoufalidis order $n$.

\end{quotation}
\vfill
\end{titlepage}

\pagebreak
\nsection{Introduction}
\label{*1}
Let $M$ be a 3-dimensional manifold with an $N$-component link $\cL$
inside it. We assign $\a_j$-dimensional irreducible representations
of $SU(2)$ to every component $\cL_j$ of $\cL$. Witten's
invariant of $M$ and $\cL$ is given~\cite{Wi}
by a path integral over all
$SU(2)$ connections $A_\mu$ on $M$:
\qq
\zaml = \int [\cD A_\mu] \exp ({ik\over 2\pi}\cs)
\pjN \Tr_{\a_j} \PexpA{\cL_j},
\label{1.1}
\qqq
here $\cs$ is the Chern-Simons action
\qq
\cs = {1\over 2} \Tr \epsilon^{\mu\nu\rho} \int_M d^3 x
\left(A_\mu\partial_\nu A_\rho + {2\over 3} A_\mu A_\nu A_\rho
\right),
\label{1.2}
\qqq
$\Tr_{\a_j} \PexpA{\cL_j}$ are traces of holonomies of $A_\mu$ along
$\cL_j$ taken in $\a_j$-dimensional representations of $SU(2)$ and
$\Tr$ of eq.~(\ref{1.2}) is the trace taken in the fundamental
2-dimensional representation. In most cases instead of the integer
number $k$ we will be using
\qq
K = k + 2.
\label{1.02}
\qqq

According to quantum field theory, the path integral~(\ref{1.1}) can
be calculated in the limit of $k\rightarrow \infty$ by the stationary
phase approximation. The stationary points of the phase~(\ref{1.2})
are flat connections. The whole path integral~(\ref{1.1}) is
presented as a sum of contributions of connected components $c$ of
the flat connection moduli space:
\qq
\zaml = \sum_c \zacl.
\label{1.10}
\qqq
Each contribution $\zacl$ is proportional to the classical
exponent $\exp(2\pi ik \csc)$, $\csc$ being the Chern-Simons action
of the flat connections of component $c$. The preexponential factor
is generally an asymptotic series in $k^{-1}$, or equivalently, in
$K^{-1}$.

Suppose that the manifold $M$ is a rational homology sphere (\Rhs).
Then the trivial connection forms a separate component of the moduli
space of flat connections. Therefore it produces a distinct
contribution to Witten's invariant~(\ref{1.10}). This contribution is
known~\cite{FrGo},~\cite{Je},~\cite{Ro*1} to be
of the following
form:
\qq
\ztrm = \frac{\sqrt{2} \pi}{K^{3\over 2} [\ordH]^{3\over 2}}
\exp \left( \sns \right),
\label{1.11}
\qqq
here $\ordH$ is the order of integer homology group. We assumed
that $M$ contained no links. We call $\snm$ ``perturbative
invariants'', because, according to quantum field theory expectations
they should be equal to the sums of $(n+1)$-loop connected Feynman
diagrams, studied, e.g. in papers~\cite{AxSi},~\cite{Ko}
and~\cite{Ta}. However caution is advised, because no direct
mathematically rigorous evidence supporting this relation has been
established yet (in fact, the results of~\cite{Ta} may even
contradict it).

In this paper we will study how $\ztrm$ changes under a rational
surgery on an algebraically split link (\ie a link with zero linking
numbers between its components) and on a boundary link. We will
derive simple surgery formulas for the invariants $\snm$ and show
that they are finite type invariants of Ohtsuki~\cite{Oh2} order $3n$
and of at most Garoufalidis~\cite{Ga} order n.

In Section~\ref{*2} we review the previous surgery formulas
of~\cite{Ro1} and~\cite{Ro3} as well as Reshetikhin's
formula~\cite{Ro2} for the Jones polynomial of a link. In
Section~\ref{*3} we derive the analog of Melvin-Morton
bound~\cite{MeMo} for the
power series expansion of the colored Jones
polynomial of algebraically split links and boundary links. By using
this bound we derive the surgery formulas for the perturbative
invariants $\snm$ for the case of a surgery on these classes of
links. The invariants are expressed in terms of derivatives of the
colored Jones polynomial and surgery data. We work out an explicit
expression for $\som$ and demonstrate
that it is consistent with J.~Hoste's
surgery formula~\cite{Ho} for Casson's invariant $\lcw$ if we
put~\cite{Ro1}
\qq
\som = 6 \lcw.
\label{1.03}
\qqq
We also show how to convert $\snm$ into integer valued invariants
$\sint$.

In Section~\ref{*4} we extend Ohtsuki and Garoufalidis definitions
of finite type invariants to rational homology spheres. We also
define an extra finite type invariant that we call \opr. We use the
results of Section~\ref{*3} to demonstrate that the perturbative
invariants $\snm$ are finite type invariants of Ohtsuki order $3n$,
\opr order $2n$ and of at most Garoufalidis order $n$. Finally, in
Section~\ref{*5} we speculate about the relation of our results to
Feynman diagram calculations of perturbative invariants
of~\cite{AxSi},~\cite{Ko} and~\cite{Ta} and to Ohtsuki's polynomial
invariant~\cite{Oh},~\cite{Oh1}.

\nsection{Surgery Formulas for Knots and General Links}
\label{*2}
\subsection{General Considerations}
\label{*2.1}

Let $\cL$ be an $N$-component link in a 3-manifold $M$. We assign
rational surgeries $\ppqj$ to its components. A rational surgery
$\ppq$ is presented by an $SL(2,\ZZ)$ matrix
\qq
\upq =
\left(
\begin{array}{cc}
p&r\\
q&s
\end{array}
\right),\qquad ps - qr = 1,
\label{1.4}
\qqq
whose coefficients show how to glue meridian and parallel of the
solid torus to meridian and parallel of the knot complement: the
meridian of the boundary of the solid torus is glued to
$p(\mbox{meridian})+q(\mbox{parallel})$ of the knot complement
(for more details see, \eg~\cite{Ro1}). We
denote by $M\p = \chi_\cL(M)$ a manifold constructed by performing
all surgeries on the components of $\cL$. The Reshetikhin-Turaev
formula~\cite{ReTu} relates Witten's invariants of $M$ and $M\p$:
\qq
\zmpk = e^{i\ffr} \sum_{1\leq \aa \leq K-1}
\zaml \pjN \utpqj_{\a_j 1},
\label{1.3}
\qqq
In this formula the matrices $\utpq_{\a\b}$ represent the group
$SL(2,\ZZ)$ in the $(K-1)$-dimensional space (of level $k$ affine
$SU(2)$ characters)~\cite{Je}:
\qq
\utpq_{\a\b} & = & i\frac {\sign{q}} {\sqrt{2K|q|}}
e^{- {i\pi\over 4} \Phi(\upq)}
\sum_{n=0}^{q-1} \sum_{\mu = \pm 1} \mu
\label{1.5}\\
&& \qquad \times
\exp \left[ \frac{i\pi}{2Kq} \left(
p\a^2 - 2\a (2Kn + \mu\b) + s (2Kn + \mu\b)^2 \right) \right],
\nonumber
\qqq
here $\Phi(\upq)$ is the Rademacher function
\qq
\Phi \left[
\begin{array}{cc}
p&r\\q&s
\end{array} \right]
= {p+s \over q} - 12 s\ppq
\label{1.6}
\qqq
and $s\ppq$ is the Dedekind sum
\qq
s\ppq = {1\over 4q} \sum_{j=1}^{q-1}
\cot\left( {\pi j\over q}\right)
\cot\left( {\pi pj\over q}\right).
\label{1.7}
\qqq
The phase $\ffr$ in eq.~(\ref{1.3}) is the framing correction
\qq
\ffr = {\pi\over 4} {K-2\over K}
\left[ \sjN \Phi(\upqj) - 3\sign{\Lt} \right],
\label{1.8}
\qqq
here $\Lt$ is an $N\times N$ matrix
\qq
\Lt_{ij} = \lij + \pqj\delta_{ij},
\label{1.9}
\qqq
$\lij$ is a linking matrix of $\cL$ and $\sign{\Lt}$ is the
difference between the number of positive and negative eigenvalues of
$\Lt$.

The formula~(\ref{1.3}) reflects the change in the whole
invariant~(\ref{1.10}). In~\cite{Ro1} we explained why its simple
modification should reflect the change in the trivial connection
contribution.
\begin{proposition}
\label{p1.1}
Suppose that a \Rhs $M\p$ is constructed by rational $\ppqj$
surgeries on the components of an $N$-component link $\cL$ in a \Rhs
$M$. Then the trivial connection contributions to Witten's invariants
of $M$ and $M\p$ are related by the following formula (\cf
eq.~(\ref{1.3})):
\qq
\ztrmpk = e^{i\ffr} {1\over 2^N} \spint{\a_j} \daa \ztram \pjN
\uhpqj_{\a_j 1}.
\label{1.12}
\qqq
Here the symbol $\spint{\a}$ means that we take only the stationary
phase contribution of the point $\a=0$ to the whole integral. The
matrix $\uhpq_{\a\b}$ is obtained from $\utpq_{\a\b}$ by
substituting $n=0$ instead of $\sum_{n=0}^{q-1}$ in eq.~(\ref{1.5}):
\qq
\uhpq_{\a\b} = \sqrt {2\over K|q|} \sign{q}
e^{-{i\pi\over 4}\Phi(\upq)}
\sin \left( {\pi\over K} {\a\b \over q} \right)
\exp \left( {i\pi \over 2Kq} (p\a^2 + s\b^2) \right).
\label{1.0012}
\qqq
\end{proposition}

If we define a trivial connection contribution to the framing
independent colored Jones polynomial by the formula
\qq
\jtram = \exp
\left( - {i\pi \over 2K} \sjN \ljj(\a_j^2 - 1) \right)
\frac{\ztram}{\ztrm},
\label{1.012}
\qqq
then the surgery formula~(\ref{1.12}) can be rewritten as
\qq
\ztrmpk & = & \ztrm e^{i\ffr} {1\over 2^N}
\spint{\a_j} \daa \jtram
\label{1.1012}
\\
&&\qquad\times
\pjN
\left( e^{{i\pi\over 2K} \ljj(\a_j^2 -1)} \uhpqj_{\a_j 1} \right).
\nonumber
\qqq

Another modification of this formula is especially useful for integer
surgeries, \ie when $p_j=1,\;\ljj=0$. It is obtained by shifting the
integration variables
\qq
\a_j \rightarrow \a_j \pm {1\over \pqlj}
\label{3.*2}
\qqq
and working with an even part of the shifted Jones polynomial
\qq
\jttram = {\prql\over 2^N}
\sum_{\mu_j=\pm 1} \left(\pjN\mu_j\right)
J^{\rm (tr)}_{\a_1 + {\mu_1 \over p_1 + q_1 l_{11}}, \ldots,
\a_N + {\mu_N \over p_N + q_N l_{NN}}}.
\label{3.*3}
\qqq
As a result, eq.~(\ref{1.1012}) transforms into
\qq
\ztrmpk & = & \ztrm \frac {i^N} {(2K)^{N\over 2}}
\left( \pjN \frac {\sign{q_j}}{\sqrt{|q_j|}} \right)
e^{-i\pi{3\over 4}\sign{\Lt}} e^{{i\pi\over K}\Dfr}
\nonumber\\
&&
\times \spint{\a_j} \daa
\exp \left[ {i\pi\over 2K} \sjN \fpqlj \a_j^2 \right]
\jttram,
\label{3.*6}
\qqq
here
\qq
\Dfr = {3\over 2}  \sign{\Lt} + {1 \over 2} \sjN
\left( 12 s\ppqj - \fpqlj - {1\over q_j (\pqlj)} \right).
\label{3.*06}
\qqq

A final comment to eqs.~(\ref{1.12}),~(\ref{1.1012}) and~(\ref{3.*6})
is that since the only flat connection on $S^3$ is trivial, then
\qq
\ztras \equiv \zasl
\label{1.2012}
\qqq
and $\jtrass$ is equal to the usual uncolored Jones polynomial
$V\left(\cL;e^{2\pi i\over K}\right)$ normalized in such a way that
$V\left(\mbox{empty link}; e^{2\pi i\over K}\right) = 1$ and
$V\left(\mbox{unknot}; e^{2\pi i\over K}\right) = 2\cos{\pi\over K}$.

In order to understand how eq.~(\ref{1.1012}) works we consider a
simple example of the stationary phase calculation:
\qq
I(K) = \spint{a} e^{iKf(a)} g(a,K) da,
\label{1.14}
\qqq
here the functions $f(a)$ and $g(a,K)$ have a smooth analytic
behavior at $a=K^{-1}=0$ and $f\p(0)=0$. We separate the quadratic
part of the exponent and then remove the odd part of the
preexponential factor, because it does not contribute to the
integral:
\qq
I(K) & = & e^{iKf(0)} \spint{a}
e^{{iK\over 2}f^{\prime\prime}(0)a^2}
\gev da,
\label{1.15}
\\
G(a,K) & = &
e^{iK \left( f(a) - f(0) - {1\over 2} f^{\prime\prime} a^2 \right)}
g(a,K),
\label{1.015}
\\
\gev & = & {1\over 2} \left( G(a,K) + G(-a,K) \right) =
\smz \snf \dmn a^{2m} K^{-n}.
\label{1.16}
\qqq
The non-zero coefficients $\dmn$ are depicted in Fig.~1. The
inequality
\qq
n \geq - {2\over 3} m
\label{1.016}
\qqq
(line $O$ in Fig.~1) comes from the fact that the expansion of the
exponent $f(a) - f(0) - {1\over 2} f^{\prime\prime}(0) a^2$ in powers
of $a$ starts with the qubic term, while the expansion of $g(a,K)$
has only negative powers of $K$.

Combining eqs.~(\ref{1.15})-(\ref{1.16}) we find that
\qq
I(K) = e^{iKf(0)} e^{{i\pi\over 4} \sign{f^{\prime\prime}(0)}}
\sqrt{2\pi \over K |f^{\prime\prime}(0)|} \sum_{n=0}^{\infty}
\Dn
K^{-n}.
\label{1.18}
\qqq
After being integrated with the gaussian factor
$e^{i {K\over 2} f^{\prime\prime}(0) a^2}$, the term $\dmn a^{2m}
K^{-n}$ contributes to the coefficient $\Delta_{m+n}$, so that
\qq
\Dn = \sum_{m=0}^{3n} {(2m)! \over m!}
\left( {i \over 2 f^{\prime\prime}(0)} \right)^m
d_{m, n-m}
\label{1.018}
\qqq
(the terms contributing to a given $\Dn$ are connected by dashed
lines on Fig.~1). The bound~(\ref{1.016}) on powers of
$K$ in power series
expansion of $\gev$ guarantees that only a finite
number of terms in that expansion is required to achieve a given
precision in expansion~(\ref{1.18}). This property makes the
stationary phase calculation of the integral~(\ref{1.14}) quite
effective. As we will see in Section~\ref{*4}, it also determines the
finite type nature of the invariants $\snm$.

\subsection{Knot Surgery Formula}
\label{*2.2}

Now we come back to the surgery formula~(\ref{1.1012}). The
substitution
\qq
\a_j = K a_j
\label{1.13}
\qqq
puts the factors $e^{\ipik \ljj (\a_j^2-1)} \uhpqj_{\a_j 1}$ of the
integrand~(\ref{1.1012}) in the form~(\ref{1.14}) with
$f(a_j) = {\pi \over 2} \fpqlj a_j^2$. It remains only to put the
Jones polynomial $\jtram$ into a similar form. If $\cL$ has only one
component, \ie if it is a knot $\cK$, then this is achieved by (the
first part of) the Melvin-Morton conjecture, which was proven by
D.~Bar-Natan and S.~Garoufalidis~\cite{BNGa} for $M=S^3$ (for a
simple path integral proof in case of a general \Rhs see~\cite{Ro1}).

\begin{proposition}
\label{p1.2}

Let $\cK$ be a knot in a \Rhs $M$. Then the trivial connection
contribution to its framing independent colored Jones polynomial has
the following expansion in powers of $\a$ and $K^{-1}$:
\qq
\jtrams = \a \snzi \sum_{0 \leq m \leq {n\over 2}}
\Dmn \a^{2m} \ipK^n,
\label{1.19}
\qqq
or equivalently,
\qq
\jtraks & = & aK \smzi \snzi \dmn (i\pi a)^{2m} \ipK^n,\;\;
a = {\a \over K},\;\; \dmn = D_{m, n+2m},
\label{1.20}
\\
\jttraks & = & \smzi \snzi \dtmn (i\pi a)^{2m} \ipK^n,
\qquad \dt_{m,0} = (2m+1) d_{m,0}.
\label{1.020}
\qqq

The dominant part of the expansions~(\ref{1.19}) and~(\ref{1.20}) is
related to the Alexander polynomial of $\cK$:
\qq
\pi a \snzi D_{n,2n} (i\pi a)^{2n} \equiv
\pi a \snzi d_{n,0} (i\pi a)^{2n} = \frac
{\sin \left( {\pi a \over m_2 d} \right)}
{\Amd}.
\label{1.21}
\qqq
The integer numbers are defined in~\cite{Ro1}, $m_2=d=1$ if $\cK$ is
homologically trivial.
\end{proposition}
The Alexander polynomial $\Am$ is normalized in such a way that
\qq
&
\displaystyle
\Delta_A\left(M,\mbox{unknot};e^{2\pi i a}\right)
 =  1,
&
\nonumber
\\
&
\displaystyle
\Am  =  \frac
{2i\sin(\pi a)} {\ordH \Rcomp},
&
\label{1.22}
\qqq
here $\Rcomp$ is the $U(1)$ Reidemeister-Ray-Singer torsion of the
knot complement $\kcomp$.

The coefficients $\dtmn$ lie above the line $G$ in Fig.~1. The
formula~(\ref{1.020}) demonstrates that the function $\jttraks$ is of
the form $g(a,K)$, that is, it has only zero or negative powers of
$K$ in its expansion. Therefore if a \Rhs $M\p$ is constructed by a
rational $\ppq$ surgery on a knot $\cK$ in a \Rhs $M$ then
\qq
\ztrmpk & = & \ztrm
i \sqrt{K \over 2} {\sign{q}\over \sqrt{|q|}}
e^{-i\pi {3\over 4} \sign{\ffpqlz}} e^{{i\pi\over K}\Dfr}
{1\over \pqlz}
\label{1.23}
\\
&&\qquad
\times \spint{a} da
\exp \left[ \ipiK \fpqlz a^2 \right]
\jttraks,
\nonumber
\\
\Dfr & = & {1 \over 2} \left[
12s\ppq - \fpqlz - {1\over q (\pqlz)} + 3\sign{\ffpqlz} \right].
\label{1.023}
\qqq
The integral is calculated similarly to the one in eq.~(\ref{1.14})
by integrating the terms of expansion~(\ref{1.020}) one by one with
the gaussian factor. The result can be expressed in terms of
invariants $\snm$ of eq.~(\ref{1.11}) if we recall that
$\ordHp = |\pqlz| \ordH$ and that $d_{00}=\tilde{d}_{00} = 1$:
\qq
&
\displaystyle
\snoi \snmp \ipK^n  =
\snoi \snm \ipK^n + {i\pi\over K}\Dfr +
\log \left( 1 + \snoi \Dn \ipK^n
\right),
&
\label{1.1023}\\
&
\displaystyle
\Dn  =  \sum_{m=0}^n (-1)^m
{(2m)!\over 2^m m!}
\left( {q\over \pqlz} \right)^m \dt_{m,n-m}.
&
\label{1.2023}
\qqq
Eq.~(\ref{1.1023}) implies that for individual perturbative
corrections
\qq
\snmp & = & \snm + \delta_{n1}\Dfr
\label{3.*10}
\\
&&\qquad
-
\sum_{m_1,\ldots,m_n \geq 0 \atop m_1 + 2m_2 + \cdots + nm_n = n}
(-1)^{m_1 + \cdots + m_n}
\frac {(m_1 + \cdots + m_n - 1)!} {m_1 ! \cdots m_n !}
\Delta_1^{m_1} \cdots \Delta_n^{m_n}.
\nonumber
\qqq

We checked in~\cite{Ro1} that the surgery formula for $S_1$ which
follows from these equations, is proportional to Walker's
formula~\cite{Wa} for Casson-Walker invariant $\lcw$ of rational
homology spheres. This led us to the relation~(\ref{1.03})

Let $\cK$ be a knot in $S^3$. The coefficients $\Dmn$ of the Taylor
expansion of its Jones polynomial are known to be Vassiliev (\ie,
finite type) invariants of order $n$. Therefore the coefficients
$\dmn=D_{m,n+2m}$ are Vassiliev invariants of order $n+2m$. The
coefficients $\dtmn$ can be expressed as linear combinations of the
coefficients $d_{m+l,n-2l}$, $l\geq 0$. Therefore $\dtmn$ is also
Vassiliev invariant of order $m+2n$.

Let $M$ be a \Rhs constructed by a rational $(p,q)$ surgery on $\cK$:
$M=\chi_\cK(S^3)$. The coefficients $\Dn$ and $\snm$ can be
considered as knot invariants of $\cK$. Eqs.~(\ref{1.2023})
and~(\ref{3.*10}) present $\snm$ as a linear combination of products
of the coefficients $\dt_{m_i,n_i}$ such that in each product
$\sum_i (n_i+m_i) = n$. As a result, for each product of
$\dt_{m_i,n_i}$ appearing in the expression of $\snm$
\qq
\sum_i (n_i+2m_i) \leq 2\sum_i (n_i + m_i) = 2n.
\label{iii.2}
\qqq
Thus we make the following conclusion:
\begin{proposition}[\cf part 1 of Question 3 of ~\cite{Ga}]
\label{pii1}
For a rational $(p,q)$ surgery on a knot $\cK\subset S^3$, the
coefficient $\snchi$, considered as an invariant of $\cK$, is a
Vassiliev invariant of order (at most) $2n$.
\end{proposition}

Note that for $l_{00}=0$ (which can be always achieved for a knot
$\cK\subset S^3$ by a suitable choice of its framing) the dependence
of $\Dn$ of $q$ in eq.~(\ref{1.2023}) is polynomial. We will use this
in order to extract the coefficients of the Alexander polynomial of
$\cK$ from the invariants $\snchi$ obtained by applying surgeries
with different values of $q$. Denote by $\chi_{\cK,(p,q)}(S^3)$ the
manifold constructed by applying the $(p,q)$ surgery to $\cK$.
Consider the surgeries $\pqmu$, $\mu_j=\pm 1$. Since
\qq
\smup \left( \qmu \right)^m =
\cases{
0 & for $m < n$ \cr
2^m\, m! & for $m=n$,\cr}
\label{iii.3}
\qqq
we find that
\qq
\smup \Delta_{n^\prime} ( \compchi ) =
\cases{
0 & for $n^\prime < n$ \cr
(-1)^n\,(2n+1)!\,p^{-n} D_{n,2n} & for $n\p=n$,\cr}
\label{iii.4}
\qqq
here we used the fact that $\dt_{n,0}=(2n+1)d_{n,0}=(2n+1)D_{n,2n}$.
Similarly, if we substitute eq.~(\ref{1.2023}) into
eq.~(\ref{3.*10}) and apply the alternating sum to both of its sides,
then eq.~(\ref{iii.3}) leads to the following
\begin{proposition}
\label{pii.2}
For a knot $\cK\subset S^3$ the alternating sum of invariants
$S_{n\p}$, $2\leq n\p \leq n$ over the surgeries $\pqmu$, $\mu_j=\pm
1$ on $\cK$ is given by the formula
\qq
\lefteqn{
\smup S_{n\p} (\compchi ) = 0\qquad  \mbox{for}\qquad n\p < n,
}
\nonumber\\
\lefteqn{
\smup S_{n} (\compchi ) }
\label{iii.6}\\
&=&
-\sum_{m_1,\ldots,m_n \geq 0 \atop m_1 + 2m_2 + \cdots + nm_n = n}
(-1)^{\sjn m_j} \left( \sjn m_j -1 \right)!\,
\pjn \left( { (-1)^j (2j + 1)! \over m_j!\,p^j } D_{j,2j}
\right)^{m_j}.
\nonumber
\qqq
here the coefficients
$D_{j,2j}$ of expansion~(\ref{1.19}) are Vassiliev invariants of
$\cK$ of order $2j$.
\end{proposition}
Eq.~(\ref{1.21}) demonstrates that $D_{j,2j}$ is a
``special'' Vassiliev invariant: it is expressed in terms of
derivatives of the Alexander polynomial
$\Delta_A (S^3,\cK;e^{2\pi ia} )$. Therefore the alternating sum
$\smup S_n (\compchi)$ is also a special Vassiliev invariant of order
$2n$ (\cf part 3 of Question 3 of ~\cite{Ga}).

One might use the relation
\qq
\Dfr = {1\over 2} q \qquad \mbox{for} \qquad l_{00}=0,\;p=1
\qqq
in order to obtain the alternating sum properties of $S_1$:
\qq
\smup S_1 ( \chi_{\cK, \left( 1, q + \sjn \mu_j \right) } ) =
\cases{
0 & for $n > 1$,\cr
q - 6 D_{1,2} & for n = 1. \cr}
\label{iii.10}
\qqq
Since $1 - 6D_{1,2}$ is proportional to the second derivative of the
Alexander polynomial, eq.~(\ref{iii.10}) is consistent with original
Casson's formula for $\lcw$ of a manifold produced by a $(1,q)$
surgery on a knot in $S^3$.

\subsection{General Link Surgery Formula}
\label{*2.3}

It is hard to present the colored Jones polynomial $\jtram$ of a
general link in the form of the integrand of eq.~(\ref{1.14}). The
closest thing to the Melvin-Morton conjecture for a general link is
Reshetikhin's formula which we proved in~\cite{Ro2} with the help of
Feynman diagrams.
\begin{proposition}
\label{p1.3}
Let $\cL$ be an $N$-component link in a \Rhs $M$. The trivial
connection contribution to its framing independent Jones polynomial
can be expressed as an integral over 3-dimensional vectors $\va_j$ of
the fixed length (\ie over the co-adjoint orbits associated with
$\a_j$-dimensional representations of $SU(2)$ which are assigned to
link components):
\qq
\jtram & = & \int_{|\va_j| = {\a_j \over K} }
\pvaj \exp \left( \ipiK \smtw \lm \right)
\label{1.25}
\\
&&\qquad
\times \pe,
\nonumber
\qqq
here $\lm$ and $\pml$ are homogeneous $SO(3)$-invariant polynomials
of degree $m$. In particular,
\qq
L_2\pva = 2 \sum_{0 \leq i < j \leq N} \lij \va_i \cdot \va_j,
\label{1.26}
\qqq
$\lij$ are the linking numbers of $\cL$.
\end{proposition}
The analog of the second part of the Melvin-Morton conjecture is the
relation~\cite{Ro2} between the polynomials $L_m$, $P_{m,l}$ and the
multicolored Alexander polynomial of the link $\cL$, defined as the
inverse of the Reidemeister-Ray-Singer torsion of its complement:
\qq
\AtN = \frac {1} {\ordH \RN}.
\label{1.27}
\qqq
\begin{proposition}
\label{p1.03}
The multicolored Alexander polynomial~(\ref{1.27}) can be expressed
in terms of polynomials $L_m$ and $P_{m,l}$:
\qq
\AtmdN & = &
-i e^{ -{i\pi\over 4} \sign{\mij}} (2\pi)^{N-2}
\frac { |\det M\pp|^{1\over 2}} {\pjN a_j}
\label{1.28}
\\
&&\qquad
\times
\left[ 1 + \smtw P_{m,0}(a_1\vn,\ldots,a_N\vn) \right]^{-1},
\nonumber
\qqq
here $\vn$ is a unit vector. A symmetric $2N\times 2N$ matrix $\mij$,
$1\leq i,j \leq N$, $\mu,\nu = 1,2$ comes from the quadratic form
\qq
\sum_{i,j=1}^N \sum_{\mu,\nu =1,2}
\mij \paa x_\mu^{(i)} x_\nu^{(j)},
\label{1.29}
\qqq
which is extracted from the exponent of eq.~(\ref{1.25})
\qq
\smtw \lm
\label{1.30}
\qqq
after performing a substitution
\qq
\va_j = a_j \left( \vn + \vx_j - {1\over 2} \vn (\vx_j^2) \right),
\qquad \vx_j \cdot \vn = 0.
\label{1.31}
\qqq
Here $x_\mu^{(j)}$ are coordinates of the vector $\vx_j$ which is
orthogonal to $\vn$. $M\pp$ is a $(N-1)\times (N-1)$ matrix obtained
from $\mij$ by crossing out two columns and two rows to which
diagonal elements $M_{ii,11}$ and $M_{jj,22}$ belong ($det M\pp$ does
not depend on the choice of $i$ and $j$). The numbers
$m_2^{(j)},d^{(j)}$ are defined in~\cite{Ro1}. If $\cL_j$ is
homologically trivial, then $m_2^{(j)}=d^{(j)}=1$.
\end{proposition}
The polynomials $\lm$ appeared to be related~\cite{Ro2} to Milnor's
linking numbers $\Mj$ of the link $\cL$. If the order $m$ Milnor's
linking numbers are well defined, then
\qq
\lm = \frac {(i\pi)^{m-2}} {m}
\sum_{1\leq j_1,\ldots,j_m \leq N} \Mj
\Tr (\vs\cdot\va_{j_1}) \cdots (\vs\cdot\va_{j_m}),
\label{1.032}
\qqq
here $\vs = (\sigma_1,\sigma_2,\sigma_3)$ is a 3d vector of Pauli
matrices. In particular,
\qq
L_3\pva & = &
- {2\pi \over 3} \sum_{i,j,k=1}^N \Mt\;
\va_i\cdot (\va_j\times\va_k),
\label{2.5}
\\
L_4\pva & = & {\pi^2\over 3} \sum_{i,j,k,l=1}^N
(\Mq - \lmu_{kijl})\;
(\va_i\times\va_j)\cdot(\va_k\times\va_l).
\label{2.6}
\qqq
We will also need the polynomial
\qq
P_{2,0}\pva = \sum_{i,j=1}^N p_{ij} \va_i\cdot\va_j
\label{2.7}
\qqq
for future use.

A combination of integrals over $\a_j = Ka_j$ in eq.~(\ref{1.1012})
with the integrals over directions of vectors $\va_j$ in
eq.~(\ref{1.25}) produces the following link surgery formula:
\begin{proposition}
\label{p1.4}
If a \Rhs $M\p$ is constructed by rational $\ppqj$ surgeries on an
$N$-component link $\cL$ in a \Rhs $M$, then
\qq
\ztrmpk & = & \ztrm (2K^3)^{N\over 2} \pjq
e^{-i\pi {3\over 4} \sign{\Lt}}
\label{1.32}
\\
&&\times
\exp \left[ \ipik \left( 3\sign{\Lt} +
\sjN \left( 12s\ppqj - \pqj - \ljj \right) \right) \right]
\nonumber\\
&&\times
\spint{\va_j} \pvtj
\exp \left[ \ipiK \left( \sum_{i,j=1}^N \Lt_{ij} \va_i\cdot\va_j
+ \smth \lm \right) \right]
\nonumber\\
&&\qquad\times
\left( \pjN \pjqs \right) \pe.
\nonumber
\qqq
\end{proposition}

By switching from integration over $a_j$ (Cartan subalgebra) to
$\va_j$ (the whole Lie algebra) we managed to put the surgery formula
in the recognizable stationary phase form~(\ref{1.14}). However we
paid a heavy price: the invariants $\snmp$ are no longer expressed in
terms of derivatives of the original colored Jones polynomial
$\jtram$, as was the case for eqs.~(\ref{1.1023}),~(\ref{1.2023}).
Instead we first have to present the polynomial in the
form~(\ref{1.25}) in order to use the coefficients of polynomials
$L_m$, $P_{m,l}$ in actual computation of the integral~(\ref{1.23}).
This is a big disadvantage of eq.~(\ref{1.23}) since we do not know
of any effective way to find the polynomials $L_m$, $P_{m,l}$ of a
link (most of them are not even unambiguously defined by
eq.~(\ref{1.25}), see~\cite{Ro2} for details). Still as we will see
in Section~\ref{*3}, for some special classes of links it is possible
to go back from eq.~(\ref{1.32}) to the formula similar to
eqs.~(\ref{1.1023}),~(\ref{1.2023}).

\nsection{Special Links}
\label{*3}
Now we will concentrate on studying the Jones polynomials and surgery
formulas of some special classes of links.

\begin{definition}
\label{d3.1}
An $N$-component link $\cL$ in a \Rhs is an algebraically split link
(\asl) if linking numbers between its components are zero:
\qq
\lij = 0, \qquad 1\leq i<j \leq N.
\label{3.01}
\qqq

A link is a special algebraically split link (\sasl) if in addition
to~(\ref{3.01}) all of its triple Milnor's linking numbers are zero:
\qq
\Mt = 0, \qquad 1\leq i,j,k \leq N.
\label{3.02}
\qqq

A link is a boundary link (\bl) if one can choose Seifert surfaces
for its components in such a way that they do not intersect.
\end{definition}
Note that all Milnor's linking numbers of a boundary link are zero.

Algebraically split links in relation to Witten's invariant were
studied by H.~Murakami and T.~Ohtsuki~\cite{Mu1}-\cite{Oh1}. In
particular, they showed that any integer homology sphere can be
constructed by integer surgeries $(1,q_j)$ on an \asl in $S^3$. They
also proved that any rational homology sphere can be constructed by
rational surgeries $(p_j,1)$ on an \asl up to a connected sum of lens
spaces. This means that instead of a desired \Rhs $M$ we may end up
with a connected sum
$M\# L_{p\p_1,1}\#\cdots\#L_{p\p_n,1}$.
This
suites our purposes since $\ztrm$ has a simple behavior under
connected sum:
\qq
\ztr(M_1\# M_2;k) = \frac
{\ztr(M_1;k)\ztr(M_2;k)} {\sqrt{2\over K} \sin \left( {\pi\over K}
\right)},
\label{3.03}
\qqq
while the trivial connection contribution to Witten's invariant of a
lens space was calculated by L.~Jeffrey~\cite{Je}:
\qq
\ztr(L_{p,q};K) = \sqrt{2\over K|p|}
\sin \left( {\pi \over K|p|} \right)
e^{-{6\pi i\over K} s(q,p)}.
\label{3.04}
\qqq

\subsection{The Jones Polynomial of Special Links}
\label{*3.1}

We are going to calculate the expansion of $\jtram$ in powers of
$\a_j$ (or equivalently, in powers of $a_j=\a_j/K$) and $K^{-1}$ with
the help of Reshetikhin's formula~(\ref{1.25}). We expand the
exponential of that formula in powers of polynomials $L_m$ and then
integrate over the directions of vectors $\va_j$ according to the
following formulas:
\qq
\int_{|\va|=a}\ida a_{\mu_1}\cdots a_{\mu_{2n+1}} & = & 0,
\label{2.1}
\\
\int_{|\va|=a}\ida a_{\mu_1}\cdots a_{\mu_{2n}} & = &
\frac {a^{2n+1}} {(2n+1)!}
\sum_{s\in S_{2n}}
\delta_{\mu_{s(1)}\mu_{s(2)}}\cdots
\delta_{\mu_{s(2n-1)}\mu_{s(2n)}},
\label{2.2}
\qqq
here $S_{2n}$ is a symmetry group of $2n$ elements.

For an \asl, $L_2 = 0$. As a result, each positive power of $K$ coming
from the expansion of the exponent of eq.~(\ref{1.25}) carries with
it at least three powers of phases $a_j$. For a \sasl, $L_2 = L_3 =
0$. Therefore each power of $K$ carries at least four powers of
phases $a_j$. Finally, for a \bl all $L_m = 0$,
so the exponential is
trivial and the power series
expansion of its Jones polynomial contains
only negative powers of $K$ (apart from the overall pre-factor $K^N$).
A ``slope index'' $\sil$ defined for special links as
\qq
\sil = \cases{
$${2\over 3}$$ & for \asl\cr
$${1\over 2}$$ & for \sasl\cr
$$0         $$ & for \bl\cr}
\label{2.02}
\qqq
allows us to formulate these results as a universal formula. It is
the analog of Melvin-Morton bound for knots.
\begin{proposition}
\label{p3.1}
The trivial connection contribution to the Jones polynomial of a
special link $\cL$ in a \Rhs $M$ has the following expansion in
powers of $K^{-1}$ and $a_j=\a_j/K$:
\qq
\jtram & = &
\left( \pjN \a_j \right) \snzi \sum_{m\leq {n\over 2-\sil}}
\Dmn (\aa) \ipK^n
\nonumber\\
& = &
K^N \left( \pjN a_j \right) \smz
\sum_{n\geq -\sil m} \dmn\paa
(i\pi)^{2m+n} \Kmn,
\label{2.3}
\qqq
while for the ``shifted'' Jones polynomial defined by
eq.~(\ref{3.*3})
\qq
\jttram =  \smz
\sum_{n\geq -\sil m} \dtmn\paa
(i\pi)^{2m+n} \Kmn.
\label{3.4}
\qqq
In these formulas $\Dmn$,
$\dmn$ and $\dtmn$ are even homogeneous polynomials
of degree $2m$:
\qq
\Dmn(\aa) & = &
\sum_{m_1,\ldots,m_N \geq 0 \atop m_1+\cdots + m_N = m}
\Dpmn \a_1^{2m_1}\cdots \a_N^{2m_N},
\label{2.301}
\\
\dmn\paa & = &
\sum_{m_1,\ldots,m_N \geq 0 \atop m_1+\cdots + m_N = m}
\dpmn a_1^{2m_1}\cdots a_N^{2m_N},
\label{2.03}
\\
\dtmn\paa & = &
\sum_{m_1,\ldots,m_N \geq 0 \atop m_1+\cdots + m_N = m}
\dtpmn a_1^{2m_1}\cdots a_N^{2m_N}.
\label{2.103}
\qqq
and $\dpmn = \Dspmn$.
\end{proposition}
The bounding lines for the polynomials $\dmn$ and $\dtmn$ in Fig.~1
are $O$ for \asl, $O\p$ for \sasl and $G$ for \bl. Note that for \asl
the polynomials of critical degrees
$d_{3m,-2m}$, $\dt_{3m,-2m}$
come exclusively from the polynomial $L_3$ of eq.~(\ref{2.5}), \ie
from the triple Milnor's linking numbers $\Mt$, while the critical
degree polynomials for \sasl
$d_{2m,-m}$, $\dt_{2m,-m}$
come exclusively from the polynomial
$L_4$ of eq.~(\ref{2.6}), \ie from quartic Milnor's linking numbers
$\Mq$.

We will need the polynomials $\dt_{1,0},\dt_{2,-1}$ and $\dt_{3,-2}$
for the surgery formula for $S_1(M)$, therefore we are going to
express them explicitly in terms of the
polynomials~(\ref{2.5})-(\ref{2.7}). The polynomial $\dt_{0,1}$ comes
from the polynomial $P_{2,0}$ in the preexponential factor of
eq.~(\ref{1.25}):
\qq
\dt_{1,0} = 3 d_{1,0} = - {3\over \pi^2}
\sjN p_{jj} a_j^2.
\label{2.8}
\qqq
The polynomial $\dt_{2,-1}$ comes from averaging the linear term in
the expansion of the exponential of the 2-color part of $L_4$
\qq
-2i\pi^3 K\sum_{1\leq i<j \leq N} \Mqij
(\va_i\times\va_j)\cdot (\va_i\times\va_j)
\label{2.08}
\qqq
over the directions of $\va_i$ and $\va_j$:
\qq
\dt_{2,-1} = 9d_{2,-1} =
 12  \sum_{1\leq i<j \leq N} \Mqij a_i^2 a_j^2.
\label{2.9}
\qqq
The polynomial $\dt_{3,-2}$ comes from averaging the quadratic term
of the expansion of the exponent of $L_3$:
\qq
\dt_{3,-2} = 27 d_{3,-2} =
-12
\sum_{1\leq i<j<k \leq N} \Mts a_i^2 a_j^2 a_k^2.
\label{2.10}
\qqq

The coefficients $p_{jj},\Mqij$ and $\Mts$ also appear as derivatives
of the Alexander polynomials of 1-, 2- and 3-component sublinks of
$\cL$. To see this, we recall~\cite{Ro2} the procedure of removing a
link component from eq.~(\ref{1.25}). To remove a particular
component $\cL_i$ of $\cL$ we have to substitute
\qq
\a_i\equiv Ka_i = 1
\label{2.011}
\qqq
and integrate over the directions of $\va_i$. We do this by expanding
the exponential of eq.~(\ref{1.25}) in all the monomials of
polynomials $L_m$ which depend on $\va_i$. As a result of the
substitution~(\ref{2.11}), these monomials become of order $K^{-1}$
or less. Therefore this procedure preserves the overall structure of
eq.~(\ref{1.25}) and leads directly to the
representation~(\ref{1.25}) of the Jones polynomial of $\cL\setminus
\cL_i$. We see that the coefficients of monomials of $L_m$, which are
independent of $\va_i$, do not change. Since $L_2 = 0$, the monomials
of $P_{2,0}$ which do not contain $\va_i$, are not modified either.

Suppose that we remove all components of $\cL$ except $\cL_j$. Then
according to eq.~(\ref{1.28}), the Alexander polynomial of $\cL_j$ in
the normalization~(\ref{1.11}) has the following
power series expansion:
\qq
\Delta_A(M,\cL_j;e^{2\pi i a_j}) =
1 - a_j^2\ppjj + O(a_j^4)
= 1 + {z^2\over 4\pi^2} \ppjj + O(z^4),
\label{2.11}
\qqq
here we used the standard variable $z$ for the Alexander polynomial:
\qq
z = -2i\sin(\pi a).
\label{2.12}
\qqq
We also assumed for simplicity that link components $\cL_j$ are
homologically trivial, so that $m_2^{(j)} = d^{(j)} = 1$. If
following J.~Hoste~\cite{Ho} we denote by $\phi_1(\cL)$ a coefficient
in front of $z^{\#\cL + 1}$ in
power series expansion of the single-colored
Alexander polynomial $\Am$ ($\#\cL$ is the number of components of
$\cL$), then according to eq.~(\ref{2.11}),
\qq
\phi(\cL_j) = {1\over 4\pi^2} \ppjj.
\label{2.13}
\qqq

Suppose that we remove all components of $\cL$, except for $\cL_i$
and $\cL_j$. Since $\lij = 0$, then according to eq.~(\ref{1.28}),
the power series
expansion of the single-colored Alexander polynomial
starts with the term
\qq
\Delta_A(M,\cL_i,\cL_j;e^{2\pi ia}) =
8i\pi^3 \Mqij a^3 + O(a^5) =
\Mqij z^3 + O(z^5),
\label{2.14}
\qqq
so that
\qq
\phi_1(\cL_i,\cL_j) = \Mqij.
\label{2.15}
\qqq
Finally for a 3-component sublink $\cL_i\bigcup \cL_j \bigcup \cL_k$
of $\cL$
\qq
\Delta_A(M,\cL_i,\cL_j,\cL_k;e^{2\pi i a}) =
16\pi^4 \Mts a^4 + O(a^6) =
\Mts z^4 + O(z^6),
\label{2.16}
\qqq
and
\qq
\phi_1(\cL_i,\cL_j,\cL_k) = \Mts.
\label{2.17}
\qqq

A combination of eqs.~(\ref{2.8}),~(\ref{2.9}),~(\ref{2.10}) with
eqs.~(\ref{2.13}),~(\ref{2.15}),~(\ref{2.17}) leads to the following
relations between the derivatives of the shifted Jones polynomial and
the derivatives of the Alexander polynomials of sublinks:
\qq
\dt_{1,0} & =  & - 12
\sjN  \left( \phi_1(\cL_j) - {1\over 24} \right)a_j^2,
\label{3.105}
\\
\dt_{2,-1} & = & 12
\sum_{1\leq i<j \leq N} \phi_1(\cL_i,\cL_j) a_i^2 a_j^2,
\label{3.205}
\\
\dt_{3,-2} & = & - 12
\sum_{1\leq i<j<k\leq N} \phi_1(\cL_i,\cL_j,\cL_k) a_i^2 a_j^2 a_k^2.
\label{3.305}
\qqq

\subsection{A Surgery Formula for Special Links}
\label{*3.2}
The Proposition~\ref{p3.1} guarantees that for all three classes of
special links --- \asl, \sasl and \bl ---
the expansion of the (shifted)
Jones polynomial is similar to that of the function~(\ref{1.015}). As
a result, we can calculate the stationary phase contribution of the
point $a_j=0$ to the integral in the surgery formula
\qq
\lefteqn{
\ztrmpk  =  \ztrm\;
i^N \left( {K\over 2} \right)^{N\over 2}
\pjq
e^{-i\pi {3\over 4}\sjN \sign{\ffpqlj}}
e^{{i\pi\over K}\Dfr}
}
\label{3.6}
\\
&&\times
{1\over \prql}
\spint{a_j} \!\!
\da \exp \left[ \ipiK \sjN \fpqlj a_j^2 \right]
\jttrakN,
\nonumber
\\
\Dfr & = & {1\over 2} \sjN \left[ 12s\ppqj - \fpqlj
- {1\over q_j(\pqlj)} + 3\sign{\ffpqlj} \right]
\label{3.06}
\qqq
by substituting the expansion~(\ref{3.4}) and integrating term by
term.
\begin{proposition}
\label{3.2}
Let $M\p$ be a \Rhs obtained by $\ppqj$ rational surgeries on
components of an $N$-component special (\asl, \sasl or \bl) link
$\cL$ in a \Rhs $M$. Then the invariants $\snmp$ and $\snm$ are
related by the following equation:
\qq
&
\displaystyle
\snoi \snmp \ipK^n  =
\snoi \snm \ipK^n + {{i\pi\over K}\Dfr}
+ \log \left( 1 + \snoi \Dn \ipK^n
\right),
&
\label{3.7}\\
&
\displaystyle
\Dn  =
\sum_{m=0}^{n\over 1-\sil}
{(-1)^m\over 2^m}
\sum_{m_1,\ldots,m_N \geq 0 \atop m_1+\cdots + m_N = m}
\dt^{(m,n-m)}_{m_1,\ldots,m_N}
\pjN
{(2m_j)!\over m_j !}
\left(
{q_j\over \pqlj} \right)^{m_j} .
&
\label{3.9}
\qqq
The individual invariants are related by the formula
\qq
\snmp & = & \snm + \delta_{n1}\Dfr
\label{3.10}
\\
&&\qquad
-
\sum_{m_1,\ldots,m_n \geq 0 \atop m_1 + 2m_2 + \cdots + nm_n = n}
(-1)^{m_1 + \cdots + m_n}
\frac {(m_1 + \cdots + m_n - 1)!} {m_1 ! \cdots m_n !}
\Delta_1^{m_1} \cdots \Delta_n^{m_n}.
\nonumber
\qqq
\end{proposition}

The case of $S_1$ is especially interesting since we demonstrated
in~\cite{Ro1} that its knot surgery formula coincides with
K.~Walker's formula for Casson's invariant $\lcw$ of \Rhs if we
set eq.~(\ref{1.03}). It is obvious from eq.~(\ref{3.10}) that for
$S_1$ one needs to know only $\Delta_1$, which for a general \asl is
expressed by eq.~(\ref{3.9}) in terms of $\dt_{1,0},\dt_{2,-1}$ and
$\dt_{3,-2}$ (note that $\dt_{0,1} = d_{0,1} = 0$, this follows from
the condition $J^{\rm (tr)}_{1,\ldots,1}(M,\cL;k) = 1$). Thus
combining eqs.~(\ref{3.105})-(\ref{3.305}) with eq.~(\ref{3.9}) we
find that
\qq
S_1(M\p) & = & S_1(M) + \Dfr +12 \left(
\sjN {\phi_1(\cL_j) - {1\over 24} \over \ffpqlj}
+
\sum_{1\leq i<j \leq N}
\frac { \phi_1(\cL_i,\cL_j)}
{\left( {p_i\over q_i} + l_{ii} \right)
\left( {p_j\over q_j} + l_{jj} \right)}
\right.
\nonumber
\\
&&\qquad\qquad\qquad\qquad
\left. +
\sum_{1\leq i<j<k \leq N}
\frac { \phi_1(\cL_i,\cL_j,\cL_k)}
{\left( {p_i\over q_i} + l_{ii} \right)
\left( {p_j\over q_j} + l_{jj} \right)
\left( {p_k\over q_k} + l_{kk} \right)}
\right).
\label{3.13}
\qqq
If we recall the Dedekind sum identity
\qq
s(p,q) + s(q,p) = {p^2 + q^2 + 1 \over 12pq} - {1 \over 4} \sign{pq},
\label{3.14}
\qqq
then it is not hard to check that for an integer surgery (\ie when
$p_j=1$, $\ljj = 0$) the substitution~(\ref{1.03}) transforms
eq.~(\ref{3.13}) into J.~Hoste's formula~\cite{Ho} for Casson's
invariant
\qq
\lcw(M\p) & = & \lcw(M) + 2\left(
\sjN q_j \phi_1(\cL_j) +
\sum_{1\leq i<j \leq N} q_i q_j \phi_1(\cL_i,\cL_j)
\right.
\label{3.15}
\\
&&\qquad\left.
+
\sum_{1\leq i<j<k \leq N} q_i q_j q_k \phi_1(\cL_i,\cL_j,\cL_k)
\right).
\nonumber
\qqq
This is another confirmation of the general relation~(\ref{1.03}).

\subsection{Integer Valued Invariants}
\label{*i}
The surgery formulas~(\ref{3.7})-(\ref{3.10}) suggest that the
invariants $\snm$ are rational numbers. In fact, we can convert them
into integers by multiplying them by factors that depend only on $n$
and $\ordH$. We are going to present a rather rough estimate of the
necessary factors.

Let $M$ be a \Rhs constructed by rational $\ppqj$
surgeries on an \asl
$\cL$ in $S^3$. The framing independent colored Jones polynomial
$\jas$ has integer coefficients in front of
the positive and negative powers of $\eipk$. We will expand the
polynomial $\jas$ in powers of $K$ in two steps. First, we introduce
a variable
\qq
x = \eipk - 1 = \ipK \snzi {1\over (n+1)!} \ipK^n
\label{ii.1}
\qqq
Since
\qq
\emipk = {1\over 1+x} = \snzi (-1)^n x^n,
\label{ii.2}
\qqq
we see that both variables $\eipk$ and $\emipk$  have integer
coefficient expansions in powers of $x$. Then a simple
relation~(\ref{ii.1}) between $x$ and $\ipK$ together with
eq.~(\ref{2.3}) imply the following expansion of the Jones polynomial
\qq
\jas = \left( \pjN \a_j \right) \snzi \sum_{m\leq {3\over 4}n}
\Dtmn(\aa) x^n,
\label{ii.3}
\qqq
here $\Dtmn(\aa)$ are homogeneous polynomials of degree $2m$:
\qq
\Dtmn = \sum_{m_1,\ldots,m_N \geq 0 \atop m_1+\cdots +m_N=m}
\Dtpmn \a_1^{2m_1} \cdots \a_N^{2m_N}.
\label{ii.4}
\qqq
The polynomials
\qq
\left( \pjN \a_j \right) \sum_{m\leq {3\over 4} n} \Dtmn(\aa)
\label{ii.5}
\qqq
are odd in all their variables. They also have integer values when
all the variables $\aa$ are integer. This means that the
polynomials~(\ref{ii.5}) can be presented as sums of products of
elementary binomial polynomials of odd degree:
\qq
P_{2m+1}(\a) = \frac {\a (\a^2 - 1) (\a^2 - 4) \cdots (\a^2 - m^2)}
{(2m+1)!}.
\label{ii.6}
\qqq
In other words,
\qq
&\displaystyle
\left( \pjN \a_j \right) \sum_{m\leq {3\over 4} n} \Dtmn(\aa)=
\sum_{0\leq m_j \leq M_j \atop (1 \leq j\leq N)}
\cnm P_{2m_1 + 1}(\a_1) \cdots P_{2m_N + 1}(\a_N),
&
\label{ii.7}\\
&\displaystyle
C_{m_1,\ldots,m_N} \in \ZZ.
&
\nonumber
\qqq
Here the numbers $M_j$ are the maximum values of the powers $m_j$
appearing in eq.~(\ref{ii.4}) for all $m\leq {3\over 4} n$. Since
$\sjN m_j = m$, then the number of $m_j\neq 0$ in each term of
eq.~(\ref{ii.7}) is not greater than $m$. Therefore
$\sum_{1\leq j \leq N \atop m_j \neq 0} (2m_j + 1) \leq 3m$ and
${(3m)! \over \pjN (2m_j + 1)! } \in \ZZ$. Taking into account that
the coefficients of the polynomial $\pjN \left( (2m_j+1)!\,
P_{2m_j+1}(\a) \right)$ are integer, we conclude that
\qq
(3m)!\, \Dtpmn \in \ZZ.
\label{ii.8}
\qqq

The estimate~(\ref{ii.8}) does not mean that prime divisors of the
denominator of $\Dtpmn$ can go as high as $3m$. Indeed, the source of
denominators is $(2m+1)!$ in eq.~(\ref{ii.6}). Since $M_j\leq m$, we
see that prime divisors of the denominator of $\Dtpmn$ are less  than
$2m+1$. Another estimate can be obtained with the help of the
following
\begin{proposition}
\label{pii.1}
For the coefficients $\dpmn$ participating in power series
expansion~(\ref{2.03}) the indices $m_j$ can not be bigger than
$m+n$:
\qq
m_j \leq m+n, \qquad 0\leq j\leq N.
\label{ii.9}
\qqq
\end{proposition}
Recall that $n$ can be negative. In that case the bound~(\ref{ii.9})
is stronger than an obvious relation $m_j\leq m$.

The Proposition~\ref{pii.1} follows easily from the Proposition~2.3
of~\cite{Ro2}, which states that the power of any vector $\va_j$ in a
polynomial $\lm$ of Reshetikhin's formula~(\ref{1.25}) is not greater
than $m-2$. A calculation of the integrals in eq.~(\ref{1.25}) leads
to the inequality~(\ref{ii.9}).

Adapting the Proposition~\ref{ii.1} to the coefficients $\Dtpmn$ we
come to the following simple corollary:
\begin{corollary}
\label{cii.1}
For the coefficients $\Dtpmn$ of the polynomials $\Dtmn$,
which participate in the
expansion~(\ref{ii.3}), each index $m_j$ can not be bigger than
$n-m$:
\qq
m_j \leq n - m, \qquad 1\leq j\leq N.
\label{ii.09}
\qqq
\end{corollary}
Note that $m_j\leq m$ and therefore $m_j\leq {n\over 2}$.

A combination of this corollary with eq.~(\ref{ii.7}) leads to the
following:
\begin{proposition}
\label{piii.1}
For the coefficients $\cnm$ of eq.~(\ref{ii.7}) there is an upper
bound on the maximum value of individual indices
\qq
m_j \leq n - \sum_{i=1}^N m_i.
\label{iii.1}
\qqq
\end{proposition}

Suppose that there exists a coefficient $\cnm$ for which the
inequality~(\ref{iii.1}) is not true, say, for $m_1$. Then the
highest degree monomial of the corresponding polynomial $\cnm
P_{2m_1+1}(\a_1)\cdots P_{2m_N+1}(\a_N)$ violates the
inequality~(\ref{ii.09}). Therefore it has to be canceled by
monomials of other polynomials
$\cnmp P_{2m\p_1+1}(\a_1)\cdots P_{2m_N\p+1}(\a_N)$ for which
apparently $m\p_j \geq m_j$, $1\leq j \leq N$ and $\sjN m\p_j > \sjN
m_j$. But the index
$m_1\p$ of these monomials again violates the
inequality~(\ref{iii.1}), so we need to go to higher values of $\sjN
m_j$ for new cancellations. Since $\sjN m_j \leq n$, this process can
not be completed. The contradiction proves the proposition.

The inequality~(\ref{iii.1}) indicates that the prime divisors of
denominators of the coefficients of the polynomial in the \rhs of
eq.~(\ref{ii.7}) can not be bigger than
$2 \left(n - \sjN m_j \right) +1$.

In order to find the polynomials $\Dmn(\aa)$ of the
expansion~(\ref{2.3}) we substitute the relation~(\ref{ii.1}) between
$x$ and $K$ into eq.~(\ref{ii.3}). The contributions to the
polynomial $\Dmn$ come from the polynomials $\Dt_{m,n-l}$, $l\geq 0$:
\qq
\Dmn = \sum_{0\leq l \leq n - {4\over 3}m} C_l \Dt_{m,n-l}.
\label{ii.10}
\qqq
The numbers $C_l$ are rational, their denominators come from the
denominator $(n+1)!$ of eq.~(\ref{ii.1}). It is easy to see that
$C_l$ has a common denominator $(2l)!$. As a result,
\qq
\left[ 2n - {8\over 3}m \right]!\, (3m)!\, \Dpmn \in\ZZ.
\label{ii.11}
\qqq
here $[x]$ is the integer part of $x$.
The polynomials $\dmn$ come from $D_{m,n+2m}$:
$\dmn = D^{(m,n+2m)}_{m_1,\ldots,m_N}$. Therefore
\qq
\left[ 2n + {4\over 3}m \right] !\, (3m)!\, \dpmn \in \ZZ
\label{ii.12}
\qqq

The polynomial $\dtmn$ comes from the polynomials $d_{\lmn}$, $l\geq
0$ through the shift of eq.~(\ref{3.*3}). The coefficient $\dtpmn$
comes from the coefficients $d^{(\lmn)}_{m_1+l_1,\dots,m_N+l_N}$,
$l_j\geq 0$, $\sum_{j=1}^N l_j = l$. The bound on powers of the
series~(\ref{2.3}) implies that ${2\over 3} (m+l)\geq 2l-n$, that is,
$l\leq \left[ {3\over 4}n + {1\over 2}m \right]$ . Since
\qq
\prql = \ordH,
\label{i.3}
\qqq
we conclude that
\qq
\left[ 2n + {4\over 3}m \right] !\,
\left[ {9\over 4}n + {9\over 2}m \right]!\;
(\ordH)^{2\mmax{l_j}} \dtpmn \in \ZZ.
\label{i.4}
\qqq

Since ${2\over 3} m \geq -n$ in $\dtmn$
because of expansion~(\ref{3.4}),
eq.~(\ref{3.9}) shows that
\qq
2^{3n} (2n)!\; (9n)!\;
(\ordH)^{2\mmax{l_j} + \mmax{m_j}} \Delta_n \in
\ZZ.
\label{i.5}
\qqq
Applying the inequality~(\ref{ii.9}) to the coefficients
$d^{(m+l,n-m-2l)}_{m_1 + l_1, \ldots, m_N + l_N}$, which
produce the coefficients $\dt^{(m,n-m)}_{m_1,\ldots,m_N}$ of
eq.~(\ref{3.9}), we find that
\qq
m + l - m_j - l_j \geq m + 2l - n,
\label{i.7}
\qqq
so that
\qq
m_j + 2l_j \leq n - l + l_j \leq n.
\label{i.8}
\qqq
The equality may be achieved if $l=l_j$. Therefore
\qq
2\mmax{l_j} + \mmax{m_j} = n
\label{i.9}
\qqq
and
\qq
2^{3n} (2n)!\; (9n)!\; (\ordH)^n \Delta_n \in \ZZ.
\label{i.10}
\qqq

The smallest denominator of each coefficient of the sum of
eq.~(\ref{3.10}) divides $m_1 + \cdots + m_N$. Therefore $n!$ may be
selected as their common denominator and
\qq
2^{3n} n!\; (2n)!\; (9n)!\;
(\ordH)^n (\snm - S_n(S^3)) \in \ZZ.
\label{i.11}
\qqq
Finally, since
\qq
\snoi S_n(S^3) \ipK^n = \log \left(
{K\over \pi} \sin \left( {\pi\over K} \right) \right),
\label{i.12}
\qqq
we come to the following conclusion:
\begin{proposition}
\label{pi.1}
The modified invariants $\sint$  of an \Rhs $M$ are integer:
\qq
\sint \equiv 2^{3n} n!\; (2n)!\; (9n)!\;
(\ordH)^n \snm \in \ZZ.
\label{i.13}
\qqq
\end{proposition}

We might have been too generous in our choice of the numerical factor
$2^{3n} n!\; (2n)!\; (9n)!$
in this definition of $\sint$. However the
calculations for Seifert manifolds (see, e.g.~\cite{Ro1}) suggest
that our choice of the power of $\ordH$ is minimal.

Recall that the factor $(9n)!$ in eq.~(\ref{i.13}) originates from
$(3m)!$ in eq.~(\ref{ii.8}). We noted there that it was needed to
remove the denominators of the polynomials $P_{2m_j+1}(\a_j)$ in the
\rhs of eq.~(\ref{ii.7}), whose prime divisors did not exceed
$2 \left( n - \sjN m_j \right) +1$. The coefficients $\Delta_{n\p}$
receive contributions from the polynomials
$P_{2m_1+1}(\a_1) \cdots P_{2m_N+1}(\a_N)$ for which
$n - \sjN m_j \leq n\p$.
Therefore the factor $(9n)!$ in eq.~(\ref{i.13}) accounts for prime
divisors which are not greater than $2n+1$. Thus we see that
\qq
(\ordH)^n \snm \in \ZZ\left[ {1\over 2}, {1\over 3},\ldots, {1\over
2n+1} \right]
\label{ii.13}
\qqq
(\cf similar results for Ohtsuki's invariants~\cite{Oh},~\cite{Oh1}).

The estimate~(\ref{ii.13}) can be improved slightly if we note that
the factor $1\over 2n + 1$ comes from the highest degree polynomials
$P_{2n+1}(\a_j)$ which may appear in the \rhs of eq.~(\ref{ii.7}) for
the polynomials $\Dt_{m, n+m}$ contributing to $\snm$ through
$\Dn(M)$. In other words, the term in the \rhs of eq.~(\ref{ii.7})
containing $P_{2n+1}$ will carry a factor $\cmms$. A simple power
counting indicates that only the highest degree monomial of the
corresponding polynomial
\qq
\cmms P_{2m_1+1}(\a_1) \cdots P_{2m_N+1}(\a_N)
\label{ii.14}
\qqq
does contribute to $\snm$. Therefore we have to follow the
transformations of only the highest degree monomial
${\a_j^{2n+1}\over (2n+1)!}$ of $P_{2n+1}(\a_j)$. It moves unchanged
from
$\left( \pjN \a_j \right) \Dt_{m,n+m}$ to
$\left( \pjN \a_j \right) D_{m, n+m}$ and transforms into
${a_j^{2n+1}\over (2n+1)!}$ inside
$\left( \pjN a_j \right) d_{m, n-m}$. A transfer to $\dt_{m, n-m}$
requires a substitution of $a_j + {1\over Kp_j}$ instead of $a_j$.
The highest even power term in
${\left( a_j + {1\over Kp_j} \right)^{2n+1} \over (2n+1)!}$ is
${1\over Kp_j} {a_j^{2n}\over (2n)!}$. Thus we see that the highest
divisor of the denominator reduced to $2n$ and we can make an
improved estimate.
\begin{proposition}
\label{piii.2}
The highest divisor of denominator of
$[\ordH]^n \snm$ is $2n$:
\qq
[\ordH]^n \snm \in \ZZ \left[ {1\over 2}, {1\over 3},\ldots, {1\over
2n} \right].
\label{ii.15}
\qqq
\end{proposition}

\nsection{Finite Type Invariants}
\label{*4}
\subsection{Definitions}
\label{*4.0}

Let $\cL$ be an $N$-component link in a 3-manifold $M$. We assign
rational surgeries $\ppqj$ to all of its components. The new manifold
constructed by performing all these surgeries is denoted as
$\chi_\cL(M)$.

T.~Ohtsuki~\cite{Oh2} and S.~Garoufalidis~\cite{Ga}
gave the following
definitions of finite type invariants of integer homology spheres
(\ihs) (we add here an extra type which we call \opr):
\begin{definition}
\label{d4.1}
A topological invariant $\lambda$
of integer homology spheres is a finite type invariant of at most
Ohtsuki (\opr, Garoufalidis) order $N$ if for any
$(N+1)$-component \asl (\sasl, \bl) $\cL$ with surgeries $(\pm 1,1)$
assigned to its components in any \ihs M, the following alternating
sum over the surgeries performed on sublinks $\cL\p\subset\cL$
(including $\cL$ itself) is equal to zero:
\qq
\sum_{\cL\p \subset \cL}
(-1)^{\#\cL\p} \lambda(\chi_{\cL\p}(M)) = 0,
\label{4.1}
\qqq
here $\#\cL\p$ is the number of components of $\cL\p$.

The invariant $\lambda$ is of Ohtsuki (\opr, Garoufalidis) order $N$
($\lambda\in O_N\; (O\p_N, \;G_N)$) if $\lambda$ is of at most order
$N$ and not of at most order $N-1$.
\end{definition}

T.~Ohtsuki proved that his classes $O_1, O_2$ were empty, while
Casson's invariant of \ihs was the only representative of his class
$O_3$. S.~Garoufalidis proved that Casson's invariant was the only
representative of his class $G_1$. He also conjectured that
\qq
O_{3n+1} = O_{3n+2} = \emptyset, \qquad O_{3n} = G_n.
\label{4.01}
\qqq

We extend the Definition~\ref{d4.1} of finite type invariants to
rational homology spheres by substituting ``arbitrary rational
surgeries $\ppqj$'' instead of ``surgeries $(\pm 1, 1)$'' in that
definition. We also conjecture that
\qq
O\p_{2n+1} = \emptyset,\qquad O\p_{2n} = G_n.
\label{4.02}
\qqq

\subsection{An Upper Estimate of Finite Type Order}
\label{*4.1}

Our first goal is to show that perturbative invariants $\snm$ are of
finite type.
\begin{proposition}
\label{p4.1}
The invariants $\snm$ of a \Rhs $M$ are finite type of at most
Ohtsuki order $3n$, at most \opr order $2n$ and at most Garoufalidis
order $n$.
\end{proposition}
Our proof is based on an observation that the difference
$\snmp-\snm$,
as it comes from eqs.~(\ref{3.10}) and~(\ref{3.9}), is sensitive to
at most $\ss$ surgeries simultaneously. The word ``simultaneously''
means that this difference can be presented as a sum of terms,
each of which can be sensitive to at most $\ss$ surgeries on link
components of $\cL$.

Let $L$ be an $N$-component link in a \Rhs $M$ with rational $\ppqj$
surgeries assigned to its components. Let $\cL\p\subset \cL$ be a
sublink of $\cL$ which does not contain a particular component, say,
$\cL_1$. Consider a contribution of the term
\qq
\dz a_2^{2m_2}\cdots a_N^{2m_N},
\qquad m = m_2 + \cdots + m_N
\label{4.3}
\qqq
from the Jones polynomial of $\cL$, which does not depend on the
phase $a_1$, to the coefficients $\Delta$ of eq.~(\ref{3.7}), as
it comes from the surgery on the link $\cL\p$ or on another link
$\cL\p\bigcup\cL_1$.
\begin{lemma}
\label{l4.1}
The contribution of the term~(\ref{4.3}) to the coefficients $\Delta$
of eq.~(\ref{3.4}) is the same for a surgery on $\cL\p$ or
$\cL\p\bigcup \cL_1$.
\end{lemma}
The proof is a direct calculation with eqs.~(\ref{3.*3})
and~(\ref{3.9}). Consider a surgery on $\cL\p$. First we have to
``remove'' the link components of $\cL\setminus \cL\p$ from the Jones
polynomial of $\cL$. This is achieved by fixing the corresponding
colors: $a_j = \a_j /K = 1/K$. In particular, we set
\qq
a_1 = {1\over K}.
\label{4.4}
\qqq
This substitution transforms the term
\qq
K^N \left( \pjN a_j \right) \dz a_2^{2m_2}\cdots a_N^{2m_N}
\label{4.5}
\qqq
of eq.~(\ref{3.2}) into
\qq
K^{N-1} \left( \prod_{j=2}^N a_j \right)
\dz a_2^{2m_2}\cdots a_N^{2m_N}.
\label{4.6}
\qqq
If a surgery is performed on $\cL\p\bigcup\cL_1$ then the variable
$a_1$ is
treated differently. First we make a part of the shift~(\ref{3.*2})
and rescaling of
eq.~(\ref{3.*3}) that are related to $a_1$.
They convert the term~(\ref{4.5}) again into the
expression~(\ref{4.6}).
%
%
Then we perform an integral over $a_1$ in eq.~(\ref{3.6}) which,
according to eq.~(\ref{3.9})
has no effect at all because $m_1$=0.
This proves the lemma.

Now we count the powers. An even homogeneous polynomial $\dmn\paa$ in
power series expansion~(\ref{2.3}) is of order $2m$, so each of its
monomials~(\ref{2.03}) depends on at most $m$ different colors.
Therefore their contribution is sensitive to at most $m$ surgeries
simultaneously. As a result of the shifts~(\ref{3.*2}) of
eq.~(\ref{3.*3}), each polynomial $\dtmn$ of the power series
expansion~(\ref{3.4}) receives the contributions of the polynomials
$d_{m+l, n-2l}$, $l\geq 0$. According to eq.~(\ref{3.9}), a
polynomial $\dtmn$ contributes to $\Delta_{m+n}$. Therefore a term
$\Delta_n$ in the surgery formula~(\ref{3.7}) receives the
contributions of the polynomials $d_{m+l, n-m-2l}$, $l,m\geq 0$.
The most surgery sensitive contribution comes from the highest
possible value of $m+l$ for a given $n$. Because of the power bound
on expansion~(\ref{2.3}) it comes from the polynomial
$d_{\ss,-{\sil\over 1-\sil}n}$. Such contribution is sensitive to at
most $\ss$ surgeries simultaneously.

Now it is easy to see that the products
$\Delta_1^{m_1}\cdots\Delta_n^{m_n}$ with $m_1+2m_2+\cdots + nm_n =n$
can be presented as sums of terms, each of which is sensitive to at
most $\ss$ surgeries. Therefore if we calculate a sum
\qq
\stnlm = \sum_{\cL\p\subset\cL}(-1)^{\#\cL\p}
S_n(\chi_{\cL}\p(M))
\label{4.07}
\qqq
for an $\left(\ss + 1\right)$-component link, then each term in
eq.~(\ref{3.10}) will be insensitive to at least one surgery, so that
it will be canceled in the alternating sum. Now it only remains to
check that
\qq
\ss = \cases{
$$3n$$ & for \asl\cr
$$2n$$ & for \sasl\cr
$$n$$          & for \bl\cr}.
\label{4.107}
\qqq
This proves the Proposition~(\ref{4.1}).

\subsection{An Exact Estimate of Ohtsuki Order}
\label{*4.2}
It follows from our proof of the Proposition~\ref{p4.1} that the most
surgery sensitive contribution to an invariant $S_n(\chi_\cL(M))$
comes from the most color-diverse monomials
\qq
c^{\left(\ss, -{\sil\over 1-\sil}n \right)}_{j_1,\ldots, j_m}
a_{j_1}^2\cdots a_{j_m}^2
\label{4.2.1}
\qqq
of the polynomial $d_{\ss,-{\sil\over 1-\sil} n}$. In case of \asl
(\sasl) these monomials have a clear geometrical origin: according to
eq.~(\ref{1.25}) they come exclusively from triple (quartic) Milnor's
linking numbers. This allows us to make a precise estimate of Ohtsuki
(\opr) order of $\snm$.
\begin{proposition}
\label{p4.2.1}
The invariants $S_n$ are of Ohtsuki (\opr) order $3n$ ( $2n$ ):
\qq
S_n\in O_{3n},\qquad S_n\in O\p_{2n}.
\label{4.2.2}
\qqq
\end{proposition}
We will present the proof for Ohtsuki invariants. The proof for \opr
invariants is similar. From now on $\cL$ in a \asl and $\sil = 2/3$.

To prove the proposition we need an effective algorithm of computing
the link invariant $\stnlm$ defined by eq.~(\ref{4.07}), for the case
of $\#\cL=3n$. Since as we have observed, the only non-zero
contributions to $\stnlm$ come from triple Milnor's linking numbers of
$\cL$, we may use a simplified version of the general link surgery
formula~(\ref{1.32}) combined with eq.~(\ref{1.11})
\qq
\lefteqn{
\exp \left( \snoi \snmp \ipK^n \right)
}
\label{4.15}\\
& \approx &
\exp \left( \snoi \snm \ipK^n \right)
\left( {K\over 2} \right)^{{3\over 2}N}
\left| \pjN {\pqlj\over q_j} \right|^{3\over 2}
e^{-i\pi{3\over 4} \sjN \sign{\ffpqlj}}
\nonumber
\\
&&\times
\spint{\va_j}d^3\va_1\cdots d^3\va_N
\exp \left[ \ipiK \left(
\sjN \fpqlj \va_j^2 -
{2\over 3} \pi \sum_{i,j,k=1}^N \Mt\;
\va_i\cdot(\va_j\times\va_k)
\right) \right]
\nonumber
\qqq
Among other things we made a substitution
\qq
\pjqs \approx {\pi\over q_j}
\label{4.015}
\qqq
in eq.~(\ref{1.32}). This approximation is justified for our
purposes. It amounts to retaining only the contribution of $\dmn$
among all polynomials $d_{m+l,n-2l}$ contributing to $\dtmn$.

Taking a logarithm of the integral in eq.~(\ref{4.15}) is a standard
exercise in combinatorics of Feynman diagrams, only in this case the
combinatorics is applied to a finite dimensional
integral~(\ref{4.15}) rather than to a path integral of quantum field
theory. The difference between the old and new invariants is
presented as a sum over diagrams
\qq
\snmp - \snm = \sum_{\dgr\in\Dgrnol} \Ddgr.
\label{4.16}
\qqq
Here $\Dgrnl$ is a set of connected $n$-loop diagrams with trivalent
vertices (see Figs.~2 and~3): each vertex represents a non-zero triple
Milnor's linking number $\Mt$ while each edge represents a link
component. Since we are interested only in the contribution of the
most color-diverse monomials~(\ref{4.2.1}) coming from the
polynomials
$d_{3n,-2n}$,
we should also require
that each component of the link should be represented at most only
once as an edge in any particular diagram. As a result, each triple
Milnor's linking number will also appear at most once, except for the
diagram of Fig.~2, where the same number appears twice. Note that a
set of participating triple linking numbers completely determines the
diagram.

A contribution $\Ddgr$ of a diagram $\dgr\in \Dgrnl$ is calculated by
expanding the exponential of the qubic part of the exponent of
eq.~(\ref{4.15}) in participating vertices (we take linear terms for
all diagrams except Fig.~2) and calculating the gaussian integrals
over $d^3\va_j$. Since
%
%
\qq
\hspace{-0.35in}
\left( {K\over 2} \right)^{3\over 2}
\left| {q\over p+ql} \right| ^{-{3\over 2}}\!\!\!\!
e^{-i\pi {3\over 4}\sign{{p\over q +l}}}
\int\! d^3\va\, a_\mu a_\nu
\exp \left[ \ipiK \left( {p\over q} + l \right) \va^2 \right]
={i\over \pi K} {q\over p+ql} \delta_{\mu\nu},
\label{4.17}
\qqq
%
the contribution $\Ddgr$ of an $n$-loop diagram
is given by the formula
\qq
\Ddgr =
4^{n-1}
{1\over 1 + \delta_{n2}} W_\dgr
\left( \prod_{j\in E_\dgr} {q_j\over \pqlj} \right)
\left( \prod_{(i,j,k)\in V_\dgr} \Mt \right).
\label{4.18}
\qqq
Here $E_\dgr$ is a set of link components appearing as edges,
$V_\dgr$ is a set of triple linking numbers appearing as vertices
(with their multiplicities) and $W_\dgr$ is a group theoretical
weight factor. It is calculated by assigning antisymmetric tensors
$\epsilon_{\mu_1\mu_2\mu_3}$ (\ie Lie algebra structure constants) to
every vertex and contracting indices along the edges. This
prescription eliminates 1-particle reducible diagrams, i.e. the
diagrams that can be split in disconnected parts by removing one
edge.

Consider now the calculation of $\stnlm$ for a $3n$-component \asl
$\cL$. The $(n+1)$-loop diagrams which contribute to the difference
\qq
S_n(\chi_{\cL\p}(M)) - \snm,
\label{4.018}
\qqq
contain $3n$ edges and therefore require $3n$ link components to
saturate them. Therefore of all sublinks $\cL\p\subset\cL$ the
difference~(\ref{4.018})
is non-zero only for $\cL\p=\cL$. Then
eq.~(\ref{4.07}) implies that
\qq
\stnlm & = & (-1)^{3n} \sum_{\dgr\in\Dgrnol} \Ddgr
\label{4.20}
\\
& = &
(-4)^n
{1\over 1+\delta_{n1}}
\left(\pjN {q_j\over \pqlj} \right)
\sum_{\dgr\in \Dgrnol} W_\dgr \prod_{(i,j,k)\in V_\dgr}\Mt.
\nonumber
\qqq

Now we can prove that $S_n\in O_{3n}$. Indeed, consider an
$(n+1)$-loop diagram $\dgr$ such that $W_\dgr \neq 0$ (it is easy to
find an example). Then similarly to~\cite{Oh2} draw a $3n$-component
\asl in $S^3$ with Borromean-type junction for every vertex of
$\dgr$. Since we kept Borromean junctions to the minimum, then the
set $\Dgrnol$ of this link contains only the original diagram
$\dgr$.
Therefore the sum~(\ref{4.20}) contains only one term which is
non-zero. This proves that the invariant $S_n$ is not of at most
Ohtsuki order $3n-1$.

A similar analysis can be carried out to show that $S_n\in O\p_{2n}$.
The Feynman diagrams will have 4-valent vertices coming from the
quartic Milnor's linking numbers $\Mq$.

\nsection{Discussion}
\label{*5}

So far the only known examples of Vassiliev invariants of links
have been
the derivatives of colored Jones polynomials corresponding to various
Lie groups. Therefore one might conjecture that for rational homology
spheres the only finite type invariants will be perturbative
invariants $S_n$. Thus it is possible that the properties of $S_n$
are universal properties of finite type invariants defined by the
Definition~\ref{d4.1}. In particular, one might hope that the
relations~(\ref{4.01}),~(\ref{4.02}) which follow so naturally from
Fig.~1, are indeed true. Each dashed line in Fig.~1 represents a
finite type invariant (or, rather, a set of invariants of the same
order) of \Rhs. Its order is equal to the $m$-coordinate of the
intersection of its dashed line with Ohtsuki, \opr and Garoufalidis
boundary lines. The main reason for us to introduce the type \opr was
that similarly to $O$ and $G$, the line $O\p$ intersects all dashed
lines at integer points.

Not all Ohtsuki diagrams~\cite{Oh2} appear in our sets $\Dgrn$. We
require the diagrams to be closed (no 1-valent vertices) and
1-particle irreducible. If someone could prove that these conditions
do follow from the Definition~\ref{d4.1},
then it might be easier to
show that $G_n\subset O_{3n}$ along the lines suggested in~\cite{Ga}.

It is easy to see that our diagrams and their group weight factors
$W_\dgr$ coincide with those appearing in the Feynman diagram
calculations of~\cite{AxSi},~\cite{Ko} and~\cite{Ta}. This is what
one might expect since according to quantum field theory, the
invariants $S_n$ should come from $(n+1)$-loop Feynman diagrams. We
present here intuitive arguments on how the Feynman diagrams may
transform into the diagrams of eq.~(\ref{4.16}).

Consider a Feynman diagram, say, the one in Fig.~2, with the edges
representing the $(1,1)$ bilocal form gauge particle propagators
$\om$. The whole expression is equal to
\qq
\int d^3 x_1 d^3 x_2
\oxx \oxx \oxx.
\label{5.1}
\qqq
Suppose that we make a rational $\ppqj$ surgery on a knot $\cK_1$.
How does the propagator $\oxx$ change outside the tubular
neighborhood, on which the surgery is performed? Since $\om$ measures
the linking numbers
\qq \lk(\cK,\cK\p) = \oint_{\cK} dx_1 \oint_{\cK\p} dx_2 \oxx
\label{5.2}
\qqq
and we know how the linking
numbers change under a rational surgery on $\cK_1$
\qq
\lk(\cK, \cK\p) \rightarrow
\lk(\cK, \cK\p) +
{p\over q} \lk(\cK,\cK_1) \lk(\cK\p,\cK_1),
\label{5.3}
\qqq
we may suggest that the propagator $\oxx$ acquires an extra piece
\qq
\oxx \rightarrow \oxx +
{p\over q} \oiky{1} \oikyp{1} \oxy{1}.
\label{5.4}
\qqq
As a result, the change in the integral~(\ref{5.1}) is proportional
to
\qq
\int d^3 x_1 d^3 x_2 \oxx \oxx \oiky{1} \oikyp{1} \oxy{1}.
\label{5.5}
\qqq
The corresponding diagram is drawn in Fig.~3. It is known~\cite{BN}
to contribute to the coefficient $d_{1,0}$ of the power series
expansion~(\ref{1.20}) of $J_{Ka}(\cK_1)$ and to the second
derivative $\phi_1(\cK_1)$ of the Alexander polynomial of $\cK_1$.
These contributions are exactly cancelled by the ghost loop. Still
this relation between the diagram of Fig.~4 and the surgery change in
the diagram of Fig.~2 is in line with our expectations that the
latter represents the Casson-Walker invariant, whose surgery formula
includes, among other terms, the derivative $\phi_1(\cK_1)$.

Let us perform a second rational surgery on another knot $\cK_2$,
such
that $\lk(\cK_1, \cK_2)$ since we want to work only with \asl. Then
the change in the integral~(\ref{5.5}) comes from ``breaking'' the
second propagator $\oxx$:
\qq
\lefteqn{
\int d^3 x_1 d^3 x_2 \oxx
\left( \oiky{1} \oikyp{1} \oxy{1} \right)
}
\label{5.6}\\
&&
\hspace{2in}
\times
\left( \oiky{2} \oikyp{2} \oxy{2} \right).
\nonumber
\qqq
The corresponding diagram is drawn in Fig.~5. Note that we could not
break the propagators $\oxy{1}$ in the expression~(\ref{5.5}) because
the result would be proportional to
\qq
\oiky{1} \oiky{2} \om(y_1,y_2) = \lk(\cK_1,\cK_2).
\qqq
The diagram of Fig.~5 contributes~\cite{Ro2} to the quartic Milnor's
linking number $\lmu_{1122}$. This is what one may expect from
eq.~(\ref{3.13}) in view of relation~(\ref{2.15}).

As we make a third surgery on $\cK_3$, the change in the
integral~(\ref{5.6}) comes from breaking the last remaining
propagator $\oxx$:
\qq
\left( \int d^3x
\oiky{1} \om(x,y_1)
\oiky{2} \om(x,y_2)
\oiky{3} \om(x,y_3)
\right)^2.
\label{5.7}
\qqq
The corresponding diagram of Fig.~6 represents~\cite{Ro2} a
contribution to the triple Milnor's linking number
$\lmu_{123}$ (\cf eqs.~(\ref{2.17}) and~(\ref{3.13})).
Its square corresponds
to the original diagram of Fig.~2, but now each vertex
represents a triple Milnor's linking number rather than a cubic term
in the Chern-Simons action~(\ref{1.2}) and each edge represents a
link component instead of a propagator $\om$. Note also that after
three surgeries we ran out of propagators to break. This indicates
that the original Feynman diagram may represent the finite type
invariant of Ohtsuki order 3.

This procedure can be applied with similar results to any Feynman
diagram containing no ghost propagators $\Omega_{0,2}$. However a
complete analysis of all diagrams and all contributions (including
the interiors of tubular neighborhoods) seems to be much more
complicated. Still the careful calculations along these lines may
shed some light on the contradiction between the
relation~(\ref{1.03}) and the results of~\cite{Ta}.

Consider for a moment integer surgeries on a special link, that is,
the surgeries for which $p_j=1$, $\ljj = 0$. Then it is easy to see
that the changes in perturbative invariants described by
eqs.~(\ref{3.9}) and~(\ref{3.10}) have a polynomial dependence on the
integer numbers $q_j$. The degree of the polynomial for $S_n$ is
$\ss$. The highest degree terms come from the coefficients lying on
the dashed lines corresponding to $\Delta_n$. These facts may
lead to another definition of finite type invariants that would use
alternating sums over surgeries performed with varying values of
$q_j$.

Finally we would like to comment on the relation between perturbative
invariants $S_n$ and Ohtsuki's invariants $\lambda_n$ introduced
in~\cite{Oh} and~\cite{Oh1}.
We show in ~\cite{Ro5} that Ohtsuki's polynomial
\qq
\tau(M) = \snzi \lambda_n (q-1)^n
\label{5.8}
\qqq
is proportional to the trivial connection contribution $\ztrm$ if we
substitute $q = \exp \left( {2\pi i\over K} \right)$:
\qq
\tau(M) =
\frac { \left( {\pi\over K} \right) }
{ \sin \left( {\pi\over K} \right)}
\exp \left[ \snoi \snm \left( {i\pi\over K} \right)^n \right].
\qqq
The proof follows from the
analysis of the surgery formulas which is similar to the one
performed in this paper. The same expansion formulas~(\ref{2.3})
and~(\ref{3.4}) have to be used in conjunction with finite gaussian
sums rather than with gaussian integrals of eq.~(\ref{3.6}).

\section*{Acknowledgements}

I am very thankful to J.~Roberts for many stimulating discussions and
for attracting my attention to the results of~\cite{Mu1}
and~\cite{Mu2}. I also appreciate the numerous conversations with
D.~Freed, C.~Gordon, L.~Kauffman, N.~Reshetikhin, A.~Vaintrob,
O.~Viro, K.~Walker and their comments. I am especially thankful to
O.~Alvarez, L.~Mezincescu and R.~Nepomechie for their constant
encouragement and friendly support.

This work was supported by the National Science Foundation
under Grant No. PHY-92 09978.

\pagebreak
\section*{Figure Captions}

\begin{description}
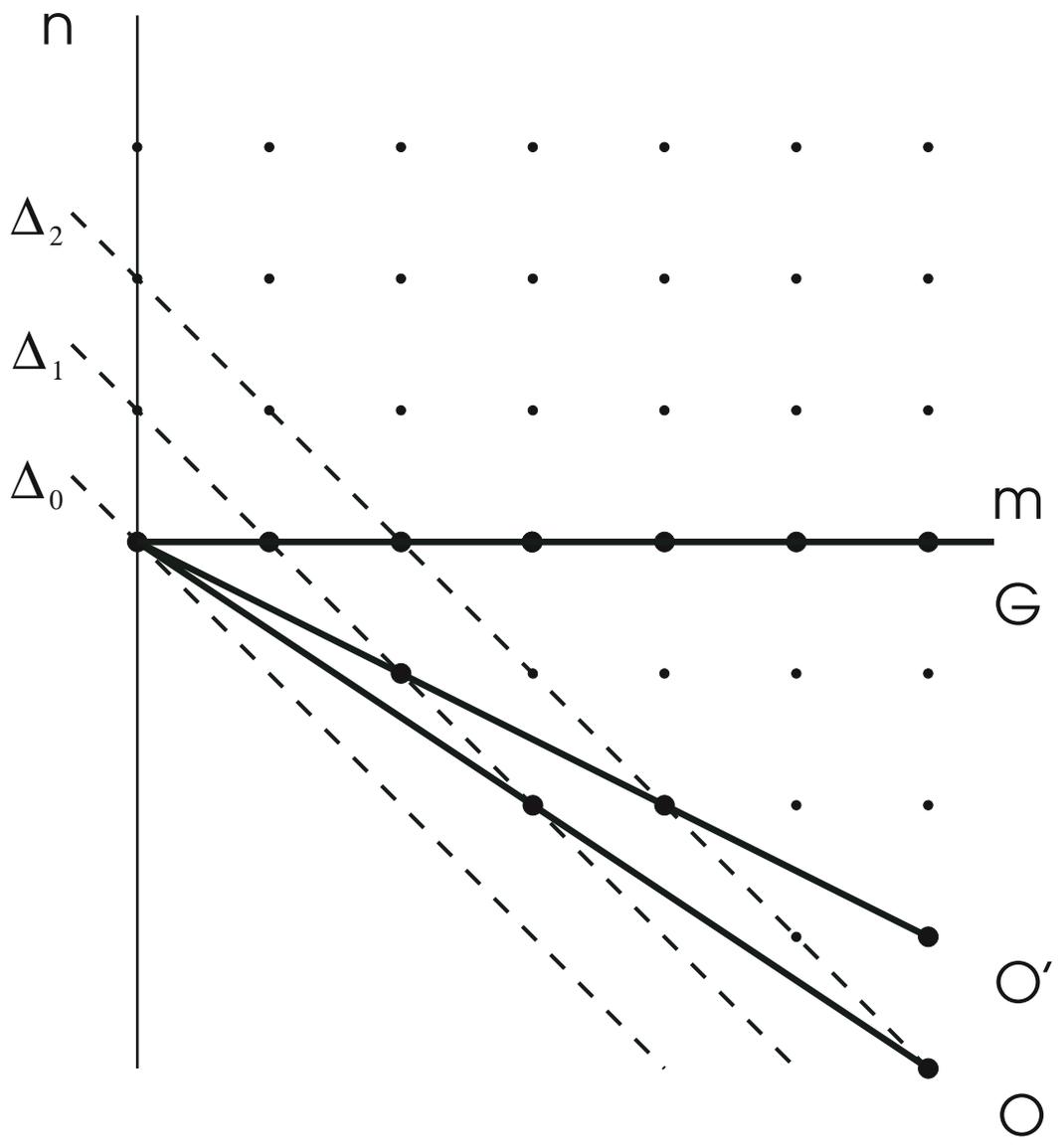
\item[Fig.~1]
The structure of power series expansions.

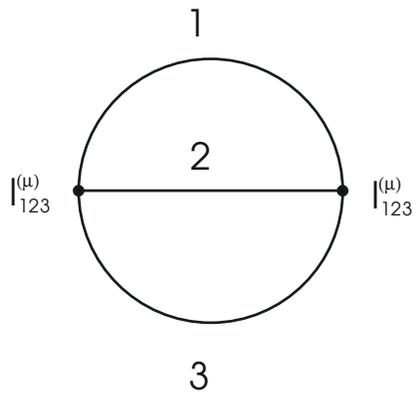
\item[Fig.~2]
The group weight diagram for $S_1$.

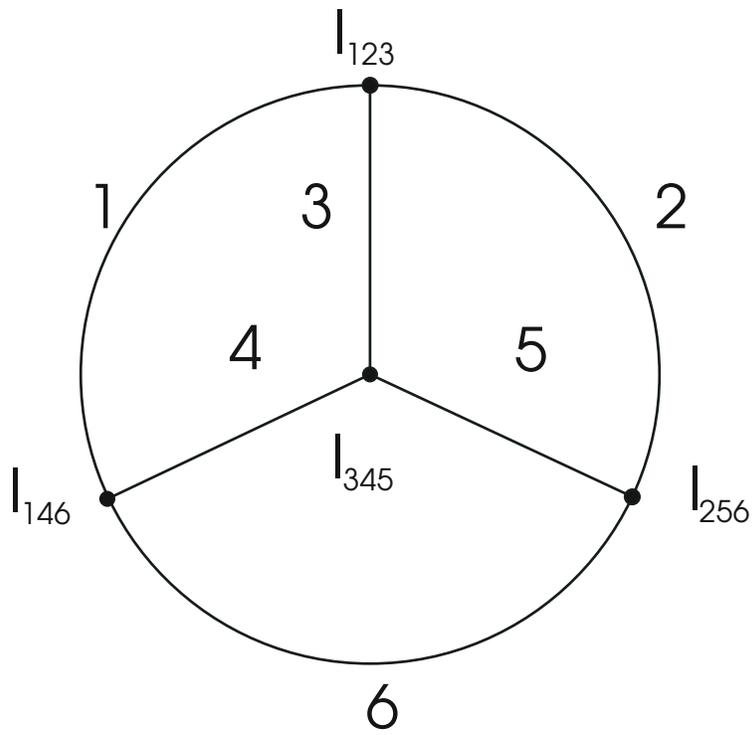
\item[Fig.~3]
A group weight diagram for $S_2$.

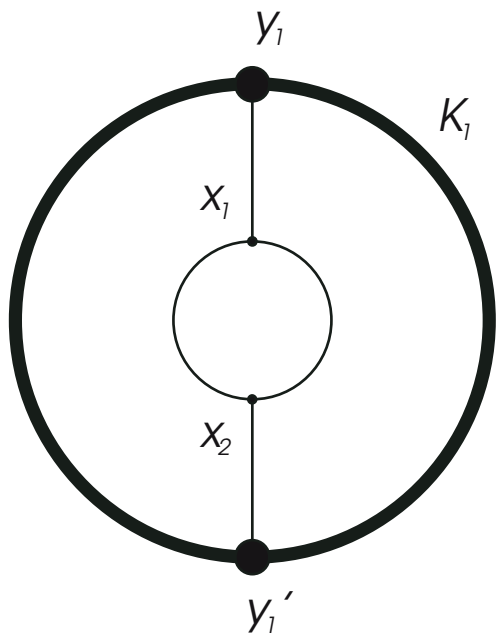
\item[Fig.~4]
A Feynman diagram contributing to $\phi_1$.

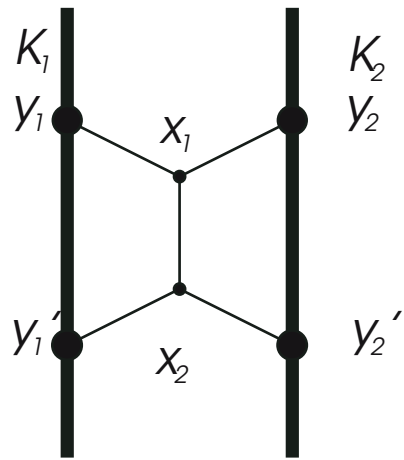
\item[Fig.~5]
A Feynman diagram contributing to $\lmu_{1122}$.

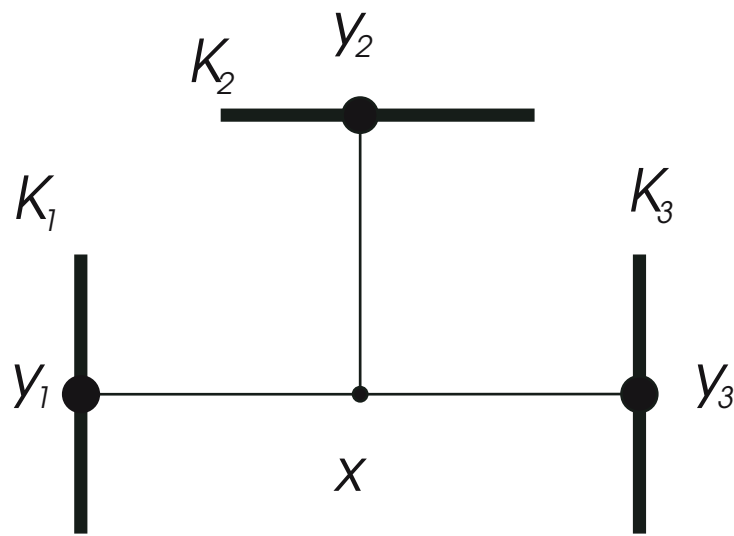
\item[Fig.~6]
A Feynman diagram contributing to $\lmu_{123}$.

\end{description}
\pagebreak

\end{document}